\documentclass[final,onefignum,onetabnum]{article}

\usepackage{array}
\usepackage[all,arc]{xy}
\usepackage{mathrsfs}
\usepackage{tcolorbox}
\usepackage{graphicx}
\usepackage{subcaption}
\usepackage{float}
\usepackage[margin=0.9 in]{geometry}
\usepackage{listings}
\usepackage{xcolor}
\usepackage{subcaption}

\usepackage[caption=false]{subfig}
\usepackage{pgfplots}
\usepackage{cite}

\usepackage{algorithmic}

\usepackage{graphicx,epstopdf}

\usepackage{subfig}
\usepackage{amsmath}
\usepackage{amssymb}
\usepackage{graphicx}
\usepackage{dcolumn}
\usepackage{mathtools}
\usepackage{amsmath,amsfonts}
\usepackage{bm}
\usepackage{amsmath}
\usepackage{amssymb}
\usepackage{color}
\usepackage{float}
\usepackage{setspace}
\usepackage{tabularx}
\usepackage{algorithm}
\usepackage{subcaption}

\allowdisplaybreaks

\newcommand{\E}{ {\mathbb{E}} }

\newcommand{\be}{\begin{equation}}
\newcommand{\ee}{\end{equation}}

\renewcommand{\d}{\mathrm{d}}

\def\ba{\begin{array}}                \def\ea{\end{array}}
\def\bel{\begin{equation}\label}      \def\ee{\end{equation}}

\colorlet{texcscolor}{blue!50!black}
\colorlet{texemcolor}{red!70!black}
\colorlet{texpreamble}{red!70!black}
\colorlet{codebackground}{black!25!white!25}

\date{}

\title{An Online Algorithm for Solving Feedback Optimal Control Problems with Partial Observations}

\author{
Siming Liang \thanks{Department of Mathematics, Florida State University, Tallahassee, Florida, USA}, Ruoyu Hu\thanks{Department of Mathematics, Florida State University, Tallahassee, Florida, USA},  Feng Bao \thanks{Department of Mathematics, Florida State University, Tallahassee, Florida, USA},  Richard Archibald \thanks{Computer Science and Mathematics Division, Oak Ridge National Laboratory, Oak Ridge, Tennessee, USA }, and Guannan Zhang \thanks{Computer Science and Mathematics Division, Oak Ridge National Laboratory, Oak Ridge, Tennessee, USA }  }

\begin{document}

\maketitle

%\tableofcontents

\begin{abstract}
This paper presents a novel methodology to tackle feedback optimal control problems in scenarios where the exact state of the controlled process is unknown. It integrates data assimilation techniques and optimal control solvers to manage partial observation of the state process, a common occurrence in practical scenarios. Traditional stochastic optimal control methods assume full state observation, which is often not feasible in real-world applications. Our approach underscores the significance of utilizing observational data to inform control policy design. Specifically, we introduce a kernel learning backward stochastic differential equation (SDE) filter to enhance data assimilation efficiency and propose a sample-wise stochastic optimization method within the stochastic maximum principle framework. Numerical experiments validate the efficacy and accuracy of our algorithm, showcasing its high efficiency in solving feedback optimal control problems with partial observation.

\end{abstract}

\section{Introduction}

In this paper, we consider a feedback optimal control problem, in which the exact state of the controlled process is not available. Stochastic optimal control is an important research subject that attracts scientists and engineers in various fields from theoretical scientific research to practical industrial production. In the stochastic optimal control problem, there is a {\it control process} (also called {\it control policy}), which controls a stochastic dynamical system whose solution called {\it state process}, and the goal of the stochastic optimal control problem is to find an optimal control process to meet some optimality conditions. For the classic stochastic optimal control problem with full observation of the state, both theoretical results and numerical methods are extensively studied. However, in practice the full observation of the state process is often not available. Instead, we have detectors/observation facilities to collect partial observational data, which provide indirect information about the state process. The theoretical formulation for partially observable stochastic optimal control is derived analytically \cite{Charalambous-98, Fleming-Pardoux-1982, Haussmann-1982, Tang-1998, Wang-Wu-Xiong-2018}, and some preliminary numerical method for the corresponding feedback control problem are developed \cite{Bao_Control_20, Bao_RL}.  To highlight the influence of data in designing control policies, in this paper we call the stochastic optimal control with partially observed controlled processes the data driven feedback control.

Since the controlled state needs to be inferred from observations, the procedure of finding the optimal control requires data analysis of observational data. The design of control actions is therefore guided by the information contained within this data. In this connection, the numerical approach for solving data-driven feedback control involves both a data assimilation procedure, which addresses the optimal filtering problem for estimating the state of the controlled dynamical system with partial noisy observations, and an optimal control solver. In this research, we present a novel approach by introducing a kernel learning backward stochastic differential equation (SDE) filter, which serves as an efficient and accurate data assimilation technique. For the stochastic optimal control problem, we propose a sample-wise stochastic optimization method within the stochastic maximum principle framework.

The optimal filtering problem for data assimilation is typically solved by recursive implementation of Bayesian inference.  The well-known approached include Kalman type filters and particle filters. The state-of-the-art Kalman fype filter is the ensemble Kalman filter (EnKF) \cite{Evense_EnKF, Tong_EnKF}. The main idea of EnKF is to use an ensemble of Kalman filter samples to characterize the probability distribution of the target state as a Gaussian distribution. As a Kalman type filter, the ensemble Kalman filter stores the information of the state variable as the mean and the covariance of Kalman filter samples. In this way, the probability density function (PDF) of the target state, which is often called the ``filtering density'', is approximated as a Gaussian distribution.  However, in nonlinear filtering problems, the filtering density is usually non-Gaussian. Therefore, the ensemble Kalman filter, which still relies on the Gaussian assumption, is not the ideal approach to solve the nonlinear filtering problem.

In addition to the ensemble Kalman filter, various effective nonlinear filtering methods have been developed to tackle nonlinearity, such as the particle filter \cite{MCMC-PF, CT1, Do2, particle-filter, Kang-PF, Sny, X_Li_Drift}, the Zakai filter \cite{Bao_zakai, zakai}, and the backward SDE filter \cite{Bao_CiCP20, BaoC20142, Bao_first,  BSDE_filter}. Among these methods, the particle filter stands out as the most widely utilized approach for addressing nonlinear filtering problems. This method utilizes a collection of Monte Carlo samples, referred to as ``particles'' to construct an empirical distribution that represents the filtering density of the target state. Upon receiving observational data, a Bayesian inference procedure is applied to assign likelihood weights to these particles, and a resampling procedure is repeatedly executed to generate additional duplicates of particles with higher weights while discarding particles with lower weights. The particle filter can handle the nonlinearity inherent in optimal filtering problems by incorporating nonlinear state dynamics into the filtering density through particle simulations. However, a significant drawback of the particle filter is the phenomenon known as ``degeneracy issue'', where particles residing in high probability regions may not adequately capture highly probable features in filtering densities. This issue becomes even more challenging when dealing with high-dimensional problems due to the ``curse of dimensionality''.  Although advanced resampling methods are proposed to address the degeneracy issue by relocating particles to high probability regions \cite{MCMC-PF, CT1, APF, Sny}, the finite particle representation for filtering density inherent in the sequential Monte Carlo framework restricts the extent of improvements that can be achieved for the particle filter in solving high-dimensional nonlinear filtering problems.

In this work, we introduce a backward SDE filter approach for solving the optimal filtering problem. Specifically, a forward backward SDE system is employed to approximate the prediction of the (controlled) state dynamics. Then, we apply Bayesian inference to incorporate the observational information and update the prediction. Numerical experiments have been conducted to show the advantageous performance of the backward SDE filter over Kalman type filters and particle filters \cite{BSDE_filter}.

\vspace{0.5em}

For solving the stochastic optimal control problem, there are two well-known approaches: the dynamic programing principle and the stochastic maximum principle \cite{Bellman1957, Feng_HJB_2013, Peng1990}.  Both approaches involve numerical simulations for large (stochastic) differential systems. In this work, we choose the stochastic maximum principle approach due to its its advantage in solving control problems with random coefficients in the state model and the fact that it could solve problems with finite dimensional terminal state constraints \cite{Yong_control}.
The mathematical foundation of the stochastic maximum principle approach is to derive a system of backward stochastic differential equation (BSDE) as the adjoint equation of the (forward) controlled state process, and we can use the solution of the adjoint BSDE to formulate the gradient of the cost function with respect to the control process. With the gradient process, the desired optimal control can be solved via gradient descent type optimization. There are several stochastic maximum principle based numerical methods to solve the stochastic optimal control problem \cite{GPM_2017, Zhao_BSDE_Control_17, Tang-1998}, and the primary computational cost in those methods lies in obtaining the numerical solution for the BSDE. Especially, when the dimension of the problem is high, a very large number of random samples are needed to describe the solution of the BSDE, which makes the numerical approximation for the gradient process very expensive \cite{Bao_AA20, ZhangJ_BSDE, Zhao_multi}. Since we aim to dynamically design control actions based on estimated parameters, efficiency  of the algorithm is very important.

To enhance efficiency, we introduce a novel approach wherein the random samples used to describe the solution of the adjoint BSDE are treated as ``data samples''. This concept underpins the development of an efficient \textit{sample-wise optimal control solver}. The essence of this solver lies in utilizing only one realization of a random sample (or a small batch of samples) to represent the solution of the BSDE and derive numerical approximations for the solution based on this single realization \cite{Bao_EAJAM20, Bao_Control_20, Bao_RL}. Additionally, the application of stochastic approximation in stochastic gradient descent optimization methodologies \cite{Convergence-SGLD, Convergence-SGD} justifies the use of sample-wise approximation for solving the BSDE in the stochastic optimal control solver \cite{Bao_SINUM_SNN}. In this way, our approach does not require solving the adjoint BSDE on a large number of random samples, shifting the computational burden from solving BSDEs to the optimization procedure for identifying the optimal control. This strategy ultimately enhances the efficiency of the stochastic optimal control solver.

 The rest of this paper is organized as follows. In Section 2, we introduce the formulation of the data driven feedback control problem and provide our backward SDE filter and stochastic maximum principle based approach. Section 3, we describe our efficient algorithm for solving the data driven feedback control in the online manner. Numerical experiments will be presented in Section 4.

\section{Problem formulation and the computational framework for the data driven feedback control problem}

\subsection{Problem setting for the data driven optimal control problem.}

In probability space $(\Omega, \mathcal{F}, \mathbb{P})$, we consider the following augmented system on time interval $[0,T]$
\begin{equation}\label{XM}%\left\{\1n\ba{ll}
\begin{aligned}
d\begin{pmatrix}S_t\\ M_t\end{pmatrix} =\begin{pmatrix}b(t,S_t,u_t )\\ g(S_t)\end{pmatrix}dt  + \begin{pmatrix}\sigma(t,S_t,u_t)&0\\0&I\end{pmatrix}
d\begin{pmatrix}W_t\\ B_t\end{pmatrix}, \qquad
\begin{pmatrix}S_0 = \xi \\ M_0 = 0 \end{pmatrix},
\end{aligned}
\end{equation}
where $S : = \{S_t\}_{t=0}^T$ is the $\mathbb{R}^d$-dimensional controlled state process with dynamics $b: [0, T] \times \mathbb{R}^d \times \mathbb{R}^m \rightarrow \mathbb{R}^d$, $\sigma: [0, T] \times \mathbb{R}^d \times \mathbb{R}^m \rightarrow \mathbb{R}^{d \times q}$ is the diffusion coefficient for a $d$-dimensional Brownian motion $W$ that perturbs the the state $S$, and $u$ is an $m$-dimensional control process valued in some set $U$ that controls the state process $S$. In the case that the state $S$ is not directly observable, we have an observation process $M$ that collects partial noisy observations on $S$ with observation function $g: \mathbb{R}^d \rightarrow \mathbb{R}^p$, and $B$ is a $p$-dimensional Brownian motion that is independent from $W$. 

Let $\mathbb{F}^B=\{\mathcal{F}_t^B\}_{t \geq 0}$ be the filtration of $B$ augmented by all the $\mathbb{P}$-null sets in $\mathcal{F}$, and $\mathbb{F}^{W,B}\equiv\{\mathcal{F}^{W,B}_t\}_{t\geq 0}$ be the filtration generated by $W$ and $B$ (augmented by $\mathbb{P}$-null sets in $\mathcal{F}$). Under mild conditions, for any square integrable random variable $\xi$ independent of $W$ and $B$, and any $\mathbb{F}^{W,B}$-progressively measurable process $u$ (valued in $U$), Eq. \eqref{XM} admits a unique solution $(S,M)$ which is $\mathbb{F}^{W,B}$-adapted. Next, we let $\mathbb{F}^M=\{\mathcal{F}^M_t\}_{t\geq 0}$ be the filtration generated by $M$ (augmented by all the $\mathbb{P}$-null sets in $\mathcal{F}$). Clearly, $\mathbb{F}^M\subset \mathbb{F}^{W,B}$, and $\mathbb{F}^M \ne \mathbb{F}^W $, $\mathbb{F}^M\ne\mathbb{F}^B$, in general. The $\mathbb{F}^M$ progressively measurable control process, denoted by $u^M$, are control actions driven by the information in observational data.

We introduce the set of data driven admissible controls as
$$\mathcal{U}_{ad}[0, T] = \left\{u^{M} : [0, T] \times \Omega  \rightarrow U \subset \mathbb{R}^m \big| u^M \text{ is } \mathbb{F}^M- \text{progressively measurable} \right\},$$
and the cost functional that measures the performance of data driven control $u^M$ is defined as
\begin{equation}\label{cost}
J(u^M) = \E\left[ \int_{0}^{T} f(t, S_t, u_t^M) dt + h(S_T) \right],
\end{equation}
where $f$ is the running cost, and $h$ is the terminal cost.

The goal of the data driven feedback control problem is to find the optimal data driven control $u^{*} \in \mathcal{U}_{ad}[0, T]$ such that
\begin{equation}
J(u^{\ast}) = \inf_{u^M \in \mathcal{U}_{ad}[0, T] } J(u^M).
\end{equation}

\subsection{Optimization for optimal control}

To solve the data driven feedback control problem, we will use the following optimization procedure derived from the stochastic maximum principle.

 When the optimal control $u^{\ast}$ is in the interior of  $\mathcal{U}_{ad}$, the gradient process of the cost functional $J^{\ast}$ with respect to the control process on time interval $t \in [0,T]$ can be derived using the Gâteaux derivative of $u^{\ast}$ and the stochastic maximum principle in the following form(see \cite{archibald2020efficient} for details):

\begin{equation}
\label{PFSGD:EQ5}
    (J^{\ast})_u^{'}(u_t^{\ast})= E\left[b_u(t,S_t^{\ast},u_t^{\ast})^{\top} \Bar{Y}_t +  \sigma_u(t,S_t^{\ast},u_t^{\ast})^{\top}\Bar{Z}_t + f_u(t,S_t^{\ast},u_t^{\ast})^{\top} | \mathcal{F}_t^M \right]
\end{equation}
stochastic processes $\Bar{Y}$ and $\Bar{Z}$ are solutions of the following forward-backward stochastic differential equations (FBSDEs) system

\begin{equation}
\label{PFSGD:EQ6}
\begin{cases}
dS_t^{\ast} = b(t,S_t^{\ast},u_t^{\ast})dt + \sigma (t,S_t^{\ast},u_t^{\ast}) dW_t, \hfill S_0 = \xi\\ 
dM_t^{\ast} = \phi(S_t)dt + dB_t, \hfill M_0 = 0  \\
d\Bar{Y}_t = (-b_x(t,S_t^{\ast},u_t^{\ast})^{\top}\Bar{Y}_t-\sigma_x(t,S_t^{\ast},u_t^{\ast})^{\top}\Bar{Z}_t-f_x(t,S_t^{\ast},u_t^{\ast})^{\top})dt\\ 
\hspace{5em}+ \Bar{Z}_t dW_t + \zeta_t dB_t,  \hfill \Bar{Y}_T=(g_x(S_T))^{\top} 
\end{cases}
\end{equation}
where $\Bar{Z}$ is the martingale representation of $\Bar{Y}$ with respect to $W$ and $\zeta$ is the martingale
representation of $\Bar{Y}$ with respect to $B$.

To solve the data driven feedback optimal control problem, we also use gradient descent type optimization and the gradient process $(J^{\ast})_u^{'}$ is defined in (\ref{PFSGD:EQ5}). Then, we can use the following gradient descent iteration to find the optimal control $u_t^{\ast}$ at any time instant $t \in [0,T]$
\begin{equation}
\label{PFSGD:EQ7}
u_t^{l+1,M} = u_t^{l,M} - r(J^{\ast})_u^{'}(u_t^{l,M}), \;\;\;\; l= 0,1,2,\dots
\end{equation}
where $r$ is the step size for the gradient. We know that the observational information $\mathcal{F}_t^M$ is gradually increased as we collect more and more data over time. Therefore, at a certain time instant $t$, we target on finding the optimal control $u_t^{\ast}$ with accessible information $\mathcal{F}_t^M$. Since the evaluation for $(J^{\ast})_u^{'}(u_t^{l,M})$ requires trajectories ${(\Bar{Y}_\tau,\Bar{Z}_\tau)}_{t \leq \tau \leq T}$ as $\Bar{Y}_t$ and $\Bar{Z}_t$ are solved backward from T to t, we take conditional expectation $E[\cdot |\mathcal{F}_t^M ]$ to the gradient process $\{(J^{\ast})_u^{'}(u_t^{l,M})\}_{t \leq \tau \leq T}$ i.e.
\begin{equation}
\label{PFSGD:EQ8}
\begin{aligned}
E[(J^{\ast})_u^{'}(u_\tau ^{l,M})|\mathcal{F}_t^M ] &= E\Big[ b_u(\tau,S_\tau,u_\tau^{l,M})^{\top}\Bar{Y}_\tau +  \sigma_u(\tau,S_\tau,u_\tau^{l,M})^{\top}\Bar{Z}_\tau \\ &+ f_u(\tau,S_\tau,u_\tau^{l,M})^{\top} |\mathcal{F}_t^M \Big], \;\;\;\;\;\;\;\; \tau \in [t,T]
\end{aligned}
\end{equation}
where $S_\tau$, $\Bar{Y}_\tau$ and $\Bar{Z}_\tau$ are corresponding to the estimated control $u_\tau^{l,M}$. for the
gradient descent iteration (\ref{PFSGD:EQ7}) on the time interval $[t, T ]$, by taking conditional expectation $E[\cdot |\mathcal{F}_t^M ]$, we obtain 
\begin{equation}
\label{PFSGD:EQ9}
E[u_\tau^{l+1,M}|\mathcal{F}_t^M ] = E[u_\tau^{l,M}|\mathcal{F}_t^M ] - rE[(J^{\ast})_u^{'}(u_\tau^{l,M})|\mathcal{F}_t^M ], \;\;\;\; l= 0,1,2,\dots \;\;\; \tau \in [t,T]
\end{equation}
When $\tau > t$, the observational information $\{\mathcal{F}_\tau^M\}_{t \leq \tau \leq T}$ is not available at time t. We use conditional expectation $E[u_\tau^{l+1,M}|\mathcal{F}_t^M ]$ to replace $u_\tau^{l,M}$ since it provides the best approximation. We denote
\begin{equation*}
    u_\tau^{l,M}|_t := E[u_\tau^{l+1,M}|\mathcal{F}_t^M ]
\end{equation*}
and then the gradient descent iteration is
\begin{equation}
\label{PFSGD:EQ10}
u_\tau^{l+1,M}|_t = u_\tau^{l,M}|_t - rE[(J^{\ast})_u^{'}(u_\tau^{l,M}|_t)|\mathcal{F}_t^M ], \;\;\;\; l= 0,1,2,\dots \;\;\; \tau \in [t,T]
\end{equation}
where $E[(J^{\ast})_u^{'}(u_\tau^{l,M}|_t)|\mathcal{F}_t^M ]$ can be obtained by solving the FBSDEs
\begin{equation}
\label{PFSGD:EQ11}
\begin{aligned}
dS_t &= b(\tau,S_\tau,u_\tau^{l,M}|_t)d\tau + \sigma (\tau,S_\tau,u_\tau^{l,M}|_t) dW_\tau, \hspace{13em} \tau \in [t,T] \\
d\Bar{Y}_t &= (-b_x(\tau,S_\tau,u_\tau^{l,M}|_t)^{\top}\Bar{Y}_\tau -\sigma_x(\tau,S_\tau,u_\tau^{l,M}|_t)^{\top}\Bar{Z}_\tau -f_x(\tau,S_\tau,u_\tau^{l,M}|_t)^{\top})d\tau \\ 
& \hspace{5em}+ \Bar{Z}_\tau dW_\tau + \zeta_\tau dB_\tau,  \hspace{15em} \Bar{Y}_T=(g_x(S_T))^{\top} 
\end{aligned}
\end{equation}
and evaluated effectively using the numerical algorithm introduced later. 

When the controlled dynamics and the observation function $\phi$ are nonlinear, we will use the optimal filtering framework to obtain the conditional expectation $E[\Psi(\tau) |\mathcal{F}_t^M ]$. For simplicity, we define 
\begin{equation}
\label{PFSGD:EQ12}
\Psi(\tau,S_\tau,u_\tau^{l,M}|_t):=b_u^{\top}(\tau,S_\tau,u_\tau^{l,M}|_t) \Bar{Y}_\tau +  \sigma_u^{\top}(\tau,S_\tau,u_\tau^{l,M}|_t)\Bar{Z}_\tau
+f_u^{\top}(\tau,S_\tau,u_\tau^{l,M}|_t)
\end{equation}
for $\tau \in [t,T]$. With the conditional pdf $p(S_t|\mathcal{F}_t^M)$ that we obtain through optimal filtering methods and the fact that $\Psi(\tau,S_\tau,u_\tau^{l,M}|_t)$ is a stochastic process depending on the state of
random variable $S_t$, the conditional gradient process $E[(J^{\ast})_u^{'}(u_t^{l,M}|_t)|\mathcal{F}_t^M ]$ in (\ref{PFSGD:EQ8}) can be obtained by the following integral
\begin{equation}
\label{PFSGD:EQ13}
E[(J^{\ast})_u^{'}(u_\tau^{l,M}|_t)|\mathcal{F}_t^M ] = \int_{\mathbb{R}^d} E[\Psi(\tau,S_\tau,u_\tau^{l,M}|_t)|S_t] \cdot p(S_t|\mathcal{F}_t^M) dx, \hspace{2em}
\tau \in [t,T]
\end{equation}
where $(\Bar{Y}_\tau, \Bar{Z}_\tau)$ in $\Psi$ can be approximated by the numerical scheme of the FBSDE system. $p(S_t|\mathcal{F}_t^M)$ can be approximated by the distribution of the backward SDE filter method introduced in the next section.

\subsection{Estimation for the state - nonlinear filtering}

In this work, we focus on the Bayes filter approach, which is carried out by recursive Bayesian estimations. In the Bayes filter, we estimate the target state on a sequence of discrete time instants $0=t_0 < t_1 < \cdots< t_{N_T}=T$ over the time interval $[0, T]$, where $N_T \in \mathbb{N}$ is the total number of time steps. The general framework of recursive Bayesian estimations is composed of two steps: a prediction step and an update step.

In the prediction step, assuming that the filtering density $p(S_{t_n} | \mathcal{F}_{t_n}^M)$ is available at the time instant $t_n$, we use the Chapman-Kolmogorov formula to propagate the dynamical model in the state equation as follows:
\begin{equation}
\label{prediction}
p(S_{t_{n+1}} | \mathcal{F}_{t_n}^M) = \int p(S_{t_{n+1}} | S_{t_n}) p(S_{t_n} | \mathcal{F}_{t_n}^M)dS_{t_n} 
\end{equation}
where $p(S_{t_{n+1}} | S_{t_n})$ is the transition probability of the state equation in (\ref{XM}), and the predicted filtering density $p(S_{t_{n+1}} | \mathcal{F}_{t_n}^M)$, which is the prior distribution in the Bayesian inference, describes the state $S$ at the time instant $t_{n+1}$ before receiving the new observational data $M_{t_{n+1}}$.

In the update step, we use the following Bayesian inference formula to incorporate the observational data into the state estimation:
\begin{equation}
\label{update}
p(S_{t_{n+1}} | \mathcal{F}_{t_{n+1}}^M) = \frac{p(M_{t_{n+1}} | S_{t_{n+1}}) p(S_{t_{n+1}} | \mathcal{F}_{t_n}^M)}{p(M_{t_{n+1}} | \mathcal{F}_{t_n}^M)}
\end{equation}
where, $p(M_{t_{n+1}} | S_{t_{n+1}})$ is the likelihood function, and the denominator $p(M_{t_{n+1}} | \mathcal{F}_{t_n}^M)$ is a normalization factor.

Then, we introduce the backward SDE filter as the theoretical preparation for our kernel learning method.

\subsubsection{The backward SDE filter}

The backward SDE filter adopts the recursive Bayesian estimations framework. The central idea is to use a system of (time-inverse) forward-backward stochastic differential equations to propagate the filtering density, and we also use Bayesian inference to incorporate the observational information into the predicted filtering density.

To proceed, we first introduce the forward-backward stochastic differential equations (FBSDEs) corresponding to the nonlinear filtering problem of the observation function in (\ref{XM}), and we consider the following FBSDEs system:
\begin{equation}
\label{BSDE:setup}
\begin{aligned}
    S_t &= S_0 + \int_0^t b(S_s) \,ds + \int_0^t \sigma_s \,dW_s, \\
    \tilde{Y}_0 &= \tilde{Y}_t - \int_0^t \tilde{Z}_s \,dW_s, \hspace{2cm} \tilde{Y}_t = \Phi(S_t),
\end{aligned}
\end{equation}
where the first equation coincides with the state equation in the nonlinear filtering problem (a standard forward SDE), and the second equation is a backward SDE. The solution of the above FBSDEs system is the pair $(\tilde{Y}, \tilde{Z})$, which is adapted to the Brownian motion $W$, i.e., $\tilde{Y}_t, \tilde{Z}_t \in \mathcal{F}_t^W$, and the solution $\tilde{Z}$ is the martingale representation of $\tilde{Y}$ with respect to $W$ \cite{peng2010backward}. $\Phi$ is a test function representing the quantity of interest. Note that the side condition of the backward SDE, i.e., $\tilde{Y}_t = \Phi(S_t)$, is given at the time instant $t$, and the solution pair $(\tilde{Y}, \tilde{Z})$ propagates backward from $t$ to $0$.

For a fixed initial state $S_0 = x \in \mathbb{R}^d$, we take the conditional expectation $E[\cdot | S_0 = x]$ on both sides of the backward SDE in (\ref{BSDE:setup}) and obtain $\tilde{Y}_0(x) = E[\Phi(S_t) | S_0 = x]$, which is a simplified version of the Feynman-Kac formula. Here, we note that the value of $\tilde{Y}_0$ is determined by the value of the state $S_0$. In addition, the solution $\tilde{Y}$ of (\ref{BSDE:setup}) is equivalent to the solution of the Kolmogorov backward equation. In other words, for the following backward parabolic type partial differential equation (PDE),

\begin{equation}
\label{BSDE:PDE}
-\frac{d u_s}{d s} = \sum_{i=1}^{d} b_i \frac{\partial u_s}{\partial x_i} + \frac{1}{2} \sum_{i,j=1}^{d} (\sigma_s \sigma_s^T )_{i,j} \frac{\partial^2 u_s}{\partial x_i \partial x_j}, \quad\quad u_t(x) = \Phi(x)
\end{equation}
where $b_i$ is the $i$-th component of the vector function $b$ in the state equation in (\ref{XM}), and we have $\tilde{Y}_0(S_0 = x) = u_0(x)$ \cite{pardoux2005backward}.

In the nonlinear filtering problem, we need to propagate the filtering density forward. The PDE that propagates the PDF of the state $S$ driven by (\ref{XM}) is the following Fokker-Planck equation:
\begin{equation}
\label{BSDE:Fokker-Planck}
\frac{dp_t}{dt} = -\sum_{i=1}^{d} \frac{\partial b_i p_t}{\partial x_i} + \frac{1}{2} \sum_{i,j=1}^{d} (\sigma_t \sigma_t^T )_{i,j} \frac{\partial^2 p_t}{\partial x_i \partial x_j} 
\end{equation}
where the initial condition $p_0$ is the distribution of the state $S_0$. Equation (\ref{BSDE:Fokker-Planck}) is the adjoint equation of the Kolmogorov backward equation, and thus, the Fokker-Planck equation is also called the Kolmogorov forward equation.

Following the analysis in \cite{pardoux2005backward} that establishes the equivalence between the FBSDEs (\ref{BSDE:setup}) and the Kolmogorov backward equation (\ref{BSDE:PDE}), one can derive that the solution $p_t$ of the above Fokker-Planck equation is equivalent to the solution $Y_t$ of the following FBSDEs:
\begin{subequations}
\label{BSDE:FBSDE}
\begin{align}
    X_0 &= X_t - \int_{0}^{t} b(X_s) \,ds + \int_{0}^{t} \sigma_s \,d\overleftarrow{W_s},  \label{BSDE:F-SDE} \\
    Y_t &= Y_0 - \int_{0}^{t} \sum_{i=1}^{d} \frac{\partial b_i}{\partial x_i}(X_s)Y_s \,ds - \int_{0}^{t} Y_s \,d\overleftarrow{W_s}, \quad Y_0 = p_0. \label{BSDE:B-SDE}
\end{align}
\end{subequations}
where the integral $\int_{0}^{t} \cdot d\overleftarrow{W_s}$ is a backward It\^o integral, integrated backward \cite{pardoux1994backward}. The equation (\ref{BSDE:F-SDE}) is a "backward SDE" as it propagates backward from $t$ to $0$, and the equation (\ref{BSDE:B-SDE}) is a "forward SDE" since its side condition is given at the time instant $0$. On the other hand, the propagation direction of (\ref{BSDE:F-SDE}) is the same as the integration direction of the backward It\^o integral. At the same time, the equation (\ref{BSDE:B-SDE}) propagates forwards given the side condition $Y_0$ = $p_0$ with a backward It\^o integral. In this way,the equations in (\ref{BSDE:FBSDE}) compose a \textit{time-inverse} FBSDEs system.

Due to the equivalence $Y_t = p_t$, the solution $Y_t$ also propagates the PDF of the state $S_t$ forward. In this connection, we use the time-inverse FBSDEs system (\ref{BSDE:FBSDE}) to predict the filtering density in the nonlinear filtering problem. Specifically, assuming we have the filtering density $p(S_{t_n} | \mathcal{F}_t^M)$ at the time step $t_n$, we solve the following time-inverse FBSDEs system:

\begin{subequations}
\label{BSDE:FBSDE2}
\begin{align}
X_{t_{n}} & =X_{t_{n+1}}-\int_{t_{n}}^{t_{n+1}} b\left(X_{s}\right) d s+\int_{t_{n}}^{t_{n+1}} \sigma_{s} d \overleftarrow{W}_{s} \label{BSDE:F-SDE2}\\
Y_{t_{n+1}}^{M_{t_{n}}} & =Y_{t_{n}}^{M_{t_{n}}}-\int_{t_{n}}^{t_{n+1}} \sum_{i=1}^{d} \frac{\partial b_{i}}{\partial x_{i}}\left(X_{s}\right) Y_{s}^{M_{t_{n}}} d s-\int_{t_{n}}^{t_{n+1}} Z_{s}^{M_{t_{n}}} d \overleftarrow{W}_{s}\label{BSDE:B-SDE2} \\
X_{t_{n+1}} & =S_{t_{n+1}}, \quad \quad \quad Y_{t_{n}}^{M_{t_{n}}}=p\left(S_{t_{n}} \mid \mathcal{F}_{t_{n}}^M\right)
\end{align}
\end{subequations}
we obtain the solution $Y_{t_{n+1}}^{M_{t_{n}}}$, which is the predicted filtering density, i.e., $Y_{t_{n+1}}^{M_{t_{n}}}=p\left(S_{t_{n+1}} \mid \mathcal{F}_{t_{n}}^M\right)$. In other words, the time-inverse FBSDEs system (\ref{BSDE:FBSDE2}) provides a mechanism to carry out the Chapman-Kolmogorov formula (\ref{prediction}) in the prediction step of the Bayes filter approach. Here, we use the superscript $M_{t_{n}}$ in $Y_{t_{n+1}}^{M_{t_{n}}}$ to emphasize that the solution $Y$ depends on the observational information $\mathcal{F}_{t_{n}}^M$.

Then, we substitute the prior distribution in the Bayesian inference (\ref{update}) by the solution $Y_{t_{n+1}}^{M_{t_{n}}}$ obtained in (\ref{BSDE:FBSDE2}) to get the posterior distribution $p\left(S_{t_{n+1}} \mid \mathcal{F}_{t_{n+1}}^M\right)$, i.e.

\begin{equation}
\label{BSDE:update}
p\left(S_{t_{n+1}} \mid \mathcal{F}_{t_{n+1}}^M\right)=\frac{p\left(M_{t_{n+1}} \mid S_{t_{n+1}}\right) Y_{t_{n+1}}^{M_{t_{n}}}}{p\left(M_{t_{n+1}} \mid \mathcal{F}_{t_{n}}^M\right)}
\end{equation}
which is the updated filtering density at time step $t_{n+1}$ that will be used for the next recursive stage.

From the above discussion, we can see that the central idea of the backward SDE filter is to use the time-inverse FBSDEs system to predict the filtering density, and then use Bayesian inference to update the predicted filtering density. As a method that carries out recursive Bayesian estimations, the backward SDE filter is also composed of a prediction step and an update step. In most practical situations, FBSDEs are not explicitly solvable. Therefore, numerical solutions for FBSDEs are needed. In the following section, we introduce a numerical algorithm to implement the above backward SDE filter framework and develop our efficient kernel learning method for the backward SDE filter.

\section{Numerical approach for data driven feedback control by BSDE-SGD}
For the numerical framework, we need the temporal partition $\Pi_{N_T}$
$$\Pi_{N_T}=\{t_n, 0=t_0 < t_1 < \cdots< t_{N_T}=T\}$$
and we use the control sequence $\{u_{t_n}^{\ast}\}_{n=1}^{N_T}$ to represent the control process $u^{\ast}$.

\subsection{Numerical schemes for the FBSDE in optimal control}

For the FBSDE system (\ref{PFSGD:EQ6}), we get the following schemes:
\begin{equation}
\begin{aligned}
\label{PF-FBSDE:EQ14}
&S_{i+1} = S_{i} + b(t_{i},S_{i},u_{t_{i}}^{l,M}|_{t_n})\triangle t_i + \sigma(t_{i},S_{i},u_{t_{i}}^{l,M}|_{t_n})\triangle W_{t_{i}}\\
&\Bar{Y}_{i} = E_i[\Bar{Y}_{i+1}] + E_i[b_x(t_{i+1},S_{i+1},u_{t_{i+1}}^{l,M}|_{t_n})^{\top}\Bar{Y}_{i+1}\\ 
& \;\;\;\;\;\;\;\;\;\;\;\;\;\;\;\;\;\;\;\;\;\;\; +\sigma_x(t_{i+1},S_{i+1},u_{t_{i+1}}^{l,M}|_{t_n})^{\top}\Bar{Z}_{i+1}+f_x(t_{i+1},S_{i+1},u_{t_{i+1}}^{l,M}|_{t_n})^{\top}] \triangle t_i\\
&\Bar{Z}_{i} =\frac{1}{\triangle t_i} E_i[\Bar{Y}_{i+1}\triangle W_{t_i}] 
\end{aligned}
\end{equation}
where $S_{i+1}, \Bar{Y}_{i}, \;and \;\Bar{Z}_{i}$ are numerical approximations for $S_{t_{i+1}}, \Bar{Y}_{t_{i}}, \;and \; \Bar{Z}_{t_{i}}$, respectively.

Then, the standard Monte Carlo method can approximate expectations with K random samples:
\begin{equation}
\begin{aligned}
\label{PF-FBSDE:EQ15}
&S_{i+1}^k = S_{i} + b(t_{i},S_{i},u_{t_{i}}^{l,M}|_{t_n})\triangle t_i + \sigma(t_{i},S_{i},u_{t_{i}}^{l,M}|_{t_n})\sqrt{\triangle t_i}\omega_i^k, k=1,2, \dots, K,\\
&\Bar{Y}_{i} =\sum_{k=1}^K \frac{\Bar{Y}_{i+1}^k}{K} + \frac{\triangle t_i}{K} \sum_{k=1}^K [b_x(t_{i+1},S_{i+1}^k,u_{t_{i+1}}^{l,M}|_{t_n})^{\top}\Bar{Y}_{i+1}^k\\ 
& \;\;\;\;\;\;\;\;\;\;\;\;\;\;\;\;\;\;\;\;\;\;\; +\sigma_x(t_{i+1},S_{i+1}^k,u_{t_{i+1}}^{l,M}|_{t_n})^{\top}\Bar{Z}_{i+1}^k+f_x(t_{i+1},S_{i+1}^k,u_{t_{i+1}}^{l,M}|_{t_n})^{\top}]\\
&\Bar{Z}_{i} =\frac{1}{\triangle t_i} \sum_{k=1}^K \frac{\Bar{Y}_{i+1}^k \sqrt{\triangle t_i}\omega_i^k}{K}
\end{aligned}
\end{equation}
where $\{ \omega_i^k \}_{k=1}^K$ is a set of random samples following the standard Gaussian distribution that
we use to describe the randomness of $\triangle W_{t_i}$.

The above schemes solve the FBSDE system (\ref{PFSGD:EQ6}) as a recursive algorithm, and the convergence of the schemes is well studied — cf. (\cite{bao2011numerical},\cite{bao2016first}, and \cite{zhao2017high}).

\subsection{Numerical algorithms for kernel learning backward SDE filter}

We first provide numerical schemes to solve the time-inverse FBSDEs system (\ref{BSDE:FBSDE2}) and then give Bayesian inference based on the numerical solution of the FBSDEs system in Subsection 3.2.1. To approximate the entire filtering density that carries the information contained in the state dynamical model and the observational data, in Subsection 3.2.2 we introduce a kernel
learning method to "learn" the filtering density from discrete density values. In Subsection 3.3 we summarize our kernel learning backward SDE filter and apply it to the stochastic optimization for the control process.

\subsubsection{Numerical schemes for time-inverse FBSDEs}
Since the equation (\ref{BSDE:F-SDE2}) is essentially an SDE with inverse propagation direction, we apply the Euler-Maruyama scheme \cite{Kloeden1992} and get
$$
X_{t_{n}}=X_{t_{n+1}}-b\left(X_{t_{n+1}}\right) \Delta t_{n}+\sigma_{t_{n+1}} \Delta W_{t_{n}}+R_{n}^{X}
$$
where $\Delta t_{n}:=t_{n+1}-t_{n}, \Delta W_{t_{n}}:=W_{t_{n+1}}-W_{t_{n}}$, and $R_{n}^{X}:=b\left(X_{t_{n+1}}\right) \Delta t_{n}-\int_{t_{n}}^{t_{n+1}} b\left(X_{s}\right) d s$ is the approximation error. By dropping the error term $R_{n}^{X}$, we obtain the following numerical scheme for (\ref{BSDE:F-SDE2}):
\begin{equation}
\label{BSDE:EulerMaruyama}
X_{n}=X_{n+1}-b\left(X_{n+1}\right) \Delta t_{n}+\sigma_{t_{n+1}} \Delta W_{t_{n}}
\end{equation}
where $X_{n}$ is the numerical approximation for $X_{t_{n}}$, and $X_{n+1}$ is a representation for the state variable $S_{t_{n+1}}$.

To solve the backward SDE (\ref{BSDE:B-SDE2}), we take conditional expectation $\mathbb{E}_{n+1}^{X}[\cdot]:=\mathbb{E}\left[\cdot \mid X_{t_{n+1}}, \mathcal{F}_{t_{n}}^M\right]$ on both sides of the equation and obtain

\begin{equation}
\label{BSDE:solveBSDE}
\mathbb{E}_{n+1}^{X}\left[Y_{t_{n+1}}^{M_{t_{n}}}\right]=\mathbb{E}_{n+1}^{X}\left[Y_{t_{n}}^{M_{t_{n}}}\right]-\int_{t_{n}}^{t_{n+1}} \mathbb{E}_{n+1}^{X}\left[\sum_{i=1}^{d} \frac{\partial b_{i}}{\partial x_{i}}\left(X_{s}\right) Y_{s}^{M_{t_{n}}}\right] d s
\end{equation}
where the backward It\^o integral $\int_{t_{n}}^{t_{n+1}} Z_{s} d \overleftarrow{W}_{s}$ is eliminated due to the martingale property of It\^o integrals. For the left hand side of the above equation, we have $Y_{t_{n+1}}^{M_{t_{n}}}=\mathbb{E}_{n+1}^{X}\left[Y_{t_{n+1}}^{M_{t_{n}}}\right]$ since $Y_{t_{n+1}}^{M_{t_{n}}}$ is adapted to $X_{t_{n+1}}$ and is $\mathcal{F}_{t_{n}}^M$ measurable. In this work, we approximate the deterministic integral on the right hand side of (\ref{BSDE:solveBSDE}) by using the right-point formula and get

\begin{equation}
\label{BSDE:solveBSDE2}
Y_{t_{n+1}}^{M_{t_{n}}}=\mathbb{E}_{n+1}^{X}\left[Y_{t_{n}}^{M_{t_{n}}}\right]-\Delta t_{n} \cdot \sum_{i=1}^{d} \frac{\partial b_{i}}{\partial x_{i}}\left(X_{t_{n+1}}\right) Y_{t_{n+1}}^{M_{t_{n}}}+R_{n}^{Y}
\end{equation}
where $R_{n}^{Y}:=\Delta t_{n} \cdot \sum_{i=1}^{d} \frac{\partial b_{i}}{\partial x_{i}}\left(X_{t_{n+1}}\right) Y_{t_{n+1}}^{M_{t_{n}}}-\int_{t_{n}}^{t_{n+1}} \mathbb{E}_{n+1}^{X}\left[\sum_{i=1}^{d} \frac{\partial b_{i}}{\partial x_{i}}\left(X_{s}\right) Y_{s}^{M_{t_{n}}}\right] d s$ is the approximation error for the integral, and we have used the fact $\sum_{i=1}^{d} \frac{\partial b_{i}}{\partial x_{i}}\left(X_{t_{n+1}}\right) Y_{t_{n+1}}^{M_{t_{n}}}=$ $\mathbb{E}_{n+1}^{X}\left[\sum_{i=1}^{d} \frac{\partial b_{i}}{\partial x_{i}}\left(X_{t_{n+1}}\right) Y_{t_{n+1}}^{M_{t_{n}}}\right]$.

Then, we drop the approximation error term $R_{n}^{Y}$ in (\ref{BSDE:solveBSDE2}) and obtain the following approxi-
mation scheme for $Y_{t_{n+1}}^{M_{t_{n}}}$ :

\begin{equation}
\label{BSDE:solveBSDE3}
Y_{n+1}^{M_{t_{n}}}=\mathbb{E}_{n+1}^{X}\left[Y_{n}^{M_{t_{n}}}\right]-\Delta t_{n} \cdot \sum_{i=1}^{d} \frac{\partial b_{i}}{\partial x_{i}}\left(X_{t_{n+1}}\right) Y_{n+1}^{M_{t_{n}}}
\end{equation}
where $Y_{n+1}^{M_{t_{n}}}$ is the approximated solution and $Y_{n}^{M_{t_{n}}}$ is an approximation of the filtering density $p\left(S_{t_{n}} \mid \mathcal{F}_{t_{n}}^M\right)$ that we obtained in the previous recursive step. We can see that the above approximation scheme is an implicit scheme. In order to calculate $Y_{n+1}^{M_{t_{n}}}$, we introduce the following fixed-point iteration procedure:

\begin{equation}
\label{BSDE:solveBSDE4}
Y_{n+1}^{M_{t_{n}}, l+1}=\mathbb{E}_{n+1}^{X}\left[Y_{n}^{M_{t_{n}}}\right]-\Delta t_{n} \cdot \sum_{i=1}^{d} \frac{\partial b_{i}}{\partial x_{i}}\left(X_{t_{n+1}}\right) Y_{n+1}^{M_{t_{n}}, l}, \quad l=0,1,2, \cdots, L-1
\end{equation}
and we let the approximated solution be $Y_{n+1}^{M_{t_{n}}}=Y_{n+1}^{M_{t_{n}}, L}$, where $L$ is a number that satisfies certain stopping criteria for the iterations, and we let the initial guess for the solution $Y_{n+1}^{M_{t_{n}}}$ be $Y_{n+1}^{M_{t_{n}}, 0}=Y_{n}^{M_{t_{n}}}$.

From the equivalence $Y_{t_{n+1}}^{M_{t_{n}}}=p\left(S_{t_{n+1}} \mid \mathcal{F}_{t_{n}}^M\right)$, the approximated solution $Y_{n+1}^{M_{t_{n}}}$ gives an approximation for the predicted filtering density. Hence the iterative scheme (\ref{BSDE:solveBSDE4}) accomplishes the prediction step in Bayesian estimation. To incorporate the observational information, we carry out the update step through Bayesian inference as follows

\begin{equation}
\label{BSDE:solveBSDE5}
\tilde{p}\left(S_{t_{n+1}} \mid \mathcal{F}_{t_{n+1}}^M\right)=\frac{p\left(M_{t_{n+1}} \mid S_{t_{n+1}}\right) Y_{n+1}^{M_{t_{n}}}}{C}
\end{equation}
where $C$ is a normalization factor, and the prior distribution is replaced by the approximation of the predicted filtering density $Y_{n+1}^{M_{t_{n}}}$. As a result, we obtain the approximated filtering density $\tilde{p}\left(S_{t_{n+1}} \mid \mathcal{F}_{t_{n+1}}^M\right)$ as desired in the backward SDE filter.

The numerical schemes (\ref{BSDE:solveBSDE4})-(\ref{BSDE:solveBSDE5}) compose a general computational framework for the backward SDE filter, which provides a recursive prediction-update mechanism that formulates the temporal propagation of the filtering density. On the other hand, the filtering density is a function that connects state positions to the probability density values at those positions. Therefore, spatial dimension approximation for the filtering density function with respect to the state variable is needed. In the following subsection, we introduce a kernel learning method to generate a continuous global approximation for the filtering density over the state space.

\subsubsection{Efficient kernel learning in the backward SDE filter}

In the scheme (\ref{BSDE:solveBSDE4}), we need to evaluate the conditional expectation $\mathbb{E}_{n+1}^{X}\left[Y_{n}^{M_{t_{n}}}\right]$ to calculate the predicted filtering density $Y_{n+1}^{M_{t_{n}}}$. So, we first discuss the approximation method for $\mathbb{E}_{n+1}^{x}\left[Y_{n}^{M_{t_{n}}}\right]$ given that $X_{t_{n+1}}=x$.

\vspace{0.5cm}
\hspace{-0.72cm} \textbf{Approximating the conditional expectation}\vspace{0.3cm}

When the target state is a high dimensional variable, Monte Carlo simulation is usually applied to evaluate high dimensional expectations. Specifically, for a given point $x \in \mathbb{R}^{d}$ in the state space, we let $X_{n+1}=x$ in the scheme (\ref{BSDE:EulerMaruyama}). Then, the conditional expectation in (\ref{BSDE:solveBSDE4}) is approximated by

\begin{equation}
\label{BSDE:ApproxCondexp}
\mathbb{E}_{n+1}^{x}\left[Y_{n}^{M_{t_{n}}}\right] \approx \hat{\mathbb{E}}_{n+1}^{x}\left[Y_{n}^{M_{t_{n}}}\right]:=\frac{\sum_{m=1}^{M} Y_{n}^{M_{t_{n}}}\left(\tilde{X}_{n}^{x, m}\right)}{M}
\end{equation}
where $M$ is the total number of Monte Carlo samples that we use to approximate the conditional expectation. The random variable $\tilde{X}_{n}^{x, m}$ in the above approximation is a Monte Carlo sample simulated by using the scheme (\ref{BSDE:EulerMaruyama}) as follows:

\begin{equation}
\label{BSDE:ApproxCondexp2}
\tilde{X}_{n}^{x, m}=x-b(x) \Delta t_{n}+\sigma_{t_{n+1}} \sqrt{\Delta t_{n}} \omega_{m}
\end{equation}
where $\omega_{m}$ is a sample drawn from the $d$-dimensional standard Gaussian distribution, and $\sqrt{\Delta t_{n}} \omega_{m}$ is the $m$-th realization of $\Delta W_{t_{n}}$. However, when the dimension of the state variable is high, the number $M$ of Monte Carlo samples needs to be very large, hence evaluating the conditional expectation $\mathbb{E}_{n+1}^{x}\left[Y_{n}^{M_{t_{n}}}\right]$ by using the Monte Carlo simulation (\ref{BSDE:ApproxCondexp}) is a computationally expensive task.

Inspired by the stochastic approximation method and its application in stochastic optimization \cite{robbins1951stochastic}, \cite{bottou2007tradeoffs}, we treat the large number of Monte Carlo samples that we use to approximate the conditional expectation $\mathbb{E}_{n+1}^{x}\left[Y_{n}^{M_{t_{n}}}\right]$ as a "large data set". Then, we adopt the methodology of stochastic approximation and use a single-sample (or a small batch of samples) to represent conditional expectations. Specifically, in each fixed-point iteration step (\ref{BSDE:solveBSDE4}), instead of using the fully-calculated Monte Carlo simulation (\ref{BSDE:ApproxCondexp}) to compute the conditional expectation, one may use one realization of the simulated sample $\tilde{X}_{n}^{x, l}$ to represent the entire set of Monte Carlo samples $\left\{\tilde{X}_{n}^{x, m}\right\}_{m=1}^{M}$ at each iteration step, where $\tilde{X}_{n}^{x, l}$ is also simulated
through (\ref{BSDE:ApproxCondexp2}) indexed by the iteration step $l$. In this way, the fixed-point iteration scheme for the approximated solution $Y_{n+1}^{M_{t_{n}}}$ at the spatial point $X_{n+1}=x$ can be carried out as follows:

\begin{equation}
\label{BSDE:ApproxCondexp3}
Y_{n+1}^{M_{t_{n}}, l+1}(x)=\tilde{\mathbb{E}}_{n+1}^{x, l}\left[Y_{n}^{M_{t_{n}}}\right]-\Delta t_{n} \cdot \sum_{i=1}^{d} \frac{\partial b_{i}}{\partial x_{i}}(x) Y_{n+1}^{M_{t_{n}}, l}, \quad l=0,1,2, \cdots, L-1
\end{equation}
where the conditional expectation $\mathbb{E}_{n+1}^{x}\left[Y_{n}^{M_{t_{n}}}\right]$ is represented by a single-sample of $Y_{n}^{M_{t_{n}}}$ corresponding to $\tilde{X}_{n}^{x, l}$, and we have

\begin{equation}
\label{BSDE:ApproxCondexp4}
\tilde{\mathbb{E}}_{n+1}^{x, l}\left[Y_{n}^{M_{t_{n}}}\right]=Y_{n}^{M_{t_{n}}}\left(\tilde{X}_{n}^{x, l}\right) .
\end{equation}

As a result, in each fixed-point iteration step, we only need to generate one sample of $X_{n}$ and evaluate the function value of the previous filtering density $Y_{n}^{M_{t_{n}}}$ at the spatial point $\tilde{X}_{n}^{x, l}$. In this way, we transfer the cost of simulating a large number of Monte Carlo samples to carrying out fixed-point iterations. Although the single-sample representation does not provide an accurate approximation for the conditional expectation, every simulated sample $\tilde{X}_{n}^{x, l}$ is effectively used to improve the estimate of the desired predicted filtering density $Y_{n+1}^{M_{t_{n}}}(x)$, which makes the overall fixed-point iteration procedure more efficient.

To use the simulated samples more effectively and to make the approximation for the conditional expectation more accurate, we modify the single-sample representation (\ref{BSDE:ApproxCondexp4}) and use the batch of samples $\left\{Y_{n}^{M_{t_{n}}}\left(\tilde{X}_{n}^{x, l}\right)\right\}_{l=0}^{L-1}$ to approximate the expectation. Precisely, we let the approximated expectation in (\ref{BSDE:ApproxCondexp3}) at each iteration step be

\begin{equation}
\label{BSDE:ApproxCondexp5}
\tilde{\mathbb{E}}_{n+1}^{x, l}\left[Y_{n}^{M_{t_{n}}}\right]=\frac{\sum_{i=1}^{l} Y_{n}^{M_{t_{n}}}\left(\tilde{X}_{n}^{x, i}\right)}{l}
\end{equation}

In this way, all the samples previously generated are used to evaluate the expectation. Note that the expectation $\mathbb{E}_{n+1}^{x}\left[Y_{n}^{M_{t_{n}}}\right]$ in the fixed-point iteration is independent of the estimation for the solution $Y_{n+1}^{M_{t_{n}}}$. Hence using more samples to approximate the expectation at each iteration step only makes the iterative scheme more accurate.

\vspace{0.2cm}
\hspace{-0.72cm} \textbf{Approximating the filtering density on random spatial sample points}\vspace{0.2cm}

By using the iterative scheme (\ref{BSDE:ApproxCondexp3}), we can calculate the predicted filtering density $Y_{n+1}^{M_{t_{n}}}$ on the state point $x$. Then, through the Bayesian inference scheme (\ref{BSDE:solveBSDE5}), we can obtain an approximation for the updated filtering density $\tilde{p}\left(S_{t_{n+1}}=x \mid \mathcal{F}_{t_{n+1}}^M\right)$ at the time step $t_{n+1}$.

To provide a complete description for the filtering density, we need to approximate the conditional PDF of the target state as a mapping from the state variable to PDF values. Standard function approximation methods use tensor-product grid points, on which we approximate function values, and then use polynomial interpolation to construct an interpolatory approximation for the entire function. If the dimension of the problem is moderately high, sparse-grid methods are often adopted as efficient alternatives to tensor-product grid interpolations. However, even advanced adaptive sparse-grid methods suffer from the "curse of dimensionality" problem. When the dimension of the problem is higher, i.e., $d \geq 10$, the cost of implementing sparse-grid approximation becomes extremely high. In many practical nonlinear filtering problems, the state dimensions are very high. Hence, applying traditional grid-based function approximation methods for filtering density is infeasible in solving high dimensional real-world problems.

An advantage of the backward SDE filter is that it allows us to approximate the filtering density on any point $x$ in the state space. In this way, we don't have to solve the problem on pre-determined meshes. Instead, in this work we use randomly generated state samples as our spatial points, and we generate spatial sample points so that they adaptively follow the conditional distribution of the target state. Specifically, assuming that we have a set of spatial points $\left\{x_{i}^{n}\right\}_{i=1}^{N}$ that follow the previous approximated filtering density $\tilde{p}\left(S_{t_{n}} \mid \mathcal{F}_{t_{n}}^M\right)$, where $N$ is the total number of spatial points, we propagate those spatial samples through the state dynamics in (\ref{XM}) by using the following Euler-Maruyama scheme:

\begin{equation}
\label{BSDE:ApproxCondexp6}
\tilde{x}_{i}^{n+1}=x_{i}^{n}+b\left(x_{i}^{n}\right) \Delta t_{n}+\sigma_{t_{n}} \sqrt{\Delta t_{n}} \omega^{i}, \quad i=1,2, \cdots, N
\end{equation}
where $\left\{\omega^{i}\right\}_{i=1}^{N}$ is a sequence of i.i.d. standard $d$-dimensional Gaussian random variables. As a result, the sample set $\left\{\tilde{x}_{i}^{n+1}\right\}_{i=1}^{N}$ forms an empirical distribution for the prior distribution $p\left(S_{t_{n+1}} \mid \mathcal{F}_{t_{n}}^M\right)$. Then, we solve the time-inverse FBSDEs (\ref{BSDE:FBSDE2}) on those spatial sample points through the scheme (\ref{BSDE:ApproxCondexp3}) to get $\left\{Y_{n+1}^{M_{t_{n}}}\left(\tilde{x}_{i}^{n+1}\right)\right\}_{i=1}^{N}$. When the new observational data $M_{t_{n+1}}$ is available, we use Bayesian inference to update the predicted filtering density and get $\left\{\tilde{p}\left(S_{t_{n+1}}=\right.\right.$ $\left.\left.\tilde{x}_{i}^{n+1} \mid \mathcal{F}_{t_{n+1}}^M\right)\right\}_{i=1}^{N}$ as our approximations for the updated filtering density. In this way, the approximations on the scattered sample points $\left\{\tilde{x}_{i}^{n+1}\right\}_{i=1}^{N}$ provide a partial description for the desired filtering density.

On the other hand, in our iterative scheme (\ref{BSDE:ApproxCondexp3}) we need the function value of $Y_{n}^{M_{t_{n}}}$ on the spatial point $\tilde{X}_{n}^{\tilde{x}_{i}^{n+1}, l}$, which is calculated from the scheme (\ref{BSDE:ApproxCondexp2}) by choosing $x=\tilde{x}_{i}^{n+1}$. Apparently, $\tilde{X}_{n}^{\tilde{x}_{i}^{n+1}, l}$ is unlikely to be one of the existing sample points, on which we have approximated
the function values for $Y_{n}^{M_{t_{n}}}$. Therefore, we need to derive an approximation for the filtering density over the entire state space. Since the filtering density is approximated on random spatial points, mesh-free methods are needed. Although traditional mesh-free interpolation methods, such as the moving least square method and the radial basis function interpolation method, could compute interpolatory approximation based on density function values at nearby spatial points, calculating the filtering density at each point separately is computationally expensive. Especially, when the dimension of the problem is high, we have to approximate the filtering density on a very large number of spatial points, which makes local approximation methods very difficult to implement.

Then, we introduce a kernel learning method to "learn" a global approximation for the entire filtering density.

\vspace{0.2cm}
\hspace{-0.72cm} \textbf{Kernel learning for the filtering density}\vspace{0.2cm}

The kernel machine utilizes the combination of a set of pre-chosen kernels to represent a model \cite{hofmann2008kernel}. For a target model in the form of a function $F$, the kernel learning method approximates the function as $F(x) \approx g\left(\sum_{k=1}^{K} \alpha_{k} \phi_{k}(x)+\beta\right), x \in \mathcal{R}^{d}$, where $\left\{\phi_{k}\right\}_{k=1}^{K}$ is a set of $K$ kernels, $g$ is an optional nonlinearity, $\left\{\alpha_{k}\right\}_{k=1}^{K}$ are weights of kernels, and $\beta$ is a bias parameter. Then we can apply the kernel learning method to approximate the filtering density under the backward SDE filter framework. Since our target function in the nonlinear filtering problem is a probability distribution, which is nonnegative and often bell-shaped, we drop the nonlinear function $g$ and the bias $\beta$ in the kernel learning model, and we choose Gaussian type functions as our kernels. Specifically, at the time step $t_{n+1}$, we use the following kernel learning scheme to formulate a global approximation for the filtering density $p\left(S_{t_{n+1}} \mid \mathcal{F}_{t_{n+1}}^M\right)$ :

\begin{equation}
\label{Kernel:density}
p_{n+1}(x):=\sum_{k=1}^{K} \alpha_{k}^{n+1} \phi_{k}^{n+1}(x), \quad x \in \mathcal{R}^{d}
\end{equation}

where $p_{n+1}$ is the kernel learned filtering density at the time step $t_{n+1}$, the Gaussian type kernel is chosen as $\phi_{k}^{n+1}(x)=\exp \left(-\left(\hat{x}_{k}^{n+1}-x\right)^{2} /\left(\lambda_{k}^{n+1}\right)^{2}\right)$ with center $\hat{x}_{k}^{n+1}$ and covariance $\lambda_{k}^{n+1}$, and the coefficient $\alpha_{k}^{n+1}>0$ is the weight of the $k$-th Gaussian kernel $\phi_{k}^{n+1}$. We can see from the scheme (\ref{Kernel:density}) that the features of the kernel learned density $p_{n+1}$ depend on the choice of kernel centers $\left\{\hat{x}_{k}^{n+1}\right\}_{k=1}^{K}$ and the parameters $\left(\alpha^{n+1}, \lambda^{n+1}\right)$, where $\alpha^{n+1}:=\left(\alpha_{1}^{n+1}, \alpha_{2}^{n+1}, \cdots, \alpha_{K}^{n+1}\right)^{T}$ and $\lambda^{n+1}:=\left(\lambda_{1}^{n+1}, \lambda_{2}^{n+1}, \cdots, \lambda_{K}^{n+1}\right)^{T}$ denote all the weights and covariances of kernels.

Since the state sample points that we generate through the scheme (\ref{BSDE:ApproxCondexp6}) adaptively follow the filtering density, which provide a good representation for the target PDF, we choose the kernel centers as a subset of the state samples $\left\{\tilde{x}_{i}^{n+1}\right\}_{i=1}^{N}$. Note that we have the value of the approximated filtering density $\tilde{p}\left(S_{t_{n+1}} \mid \mathcal{F}_{t_{n+1}}^M\right)$ on each sample point $S_{t_{n+1}}=\tilde{x}_{i}^{n+1}$. Therefore, we can choose state samples with high density values as kernel centers. To avoid using too many samples in the mode of the distribution and to capture more features of the filtering density, we use importance sampling to choose $\left\{\hat{x}_{k}^{n+1}\right\}_{k=1}^{K}$ \cite{gordon1993novel} instead of only using samples with highest density values.

To determine the parameters of the kernels in kernel learning, we use the approximated filtering density values $\left\{\tilde{p}\left(S_{t_{n+1}}=\tilde{x}_{i}^{n+1} \mid \mathcal{F}_{t_{n+1}}^M\right)\right\}_{i=1}^{N}$ as simulated "training data", and we aim to find kernel parameters so that the kernel learned filtering density $p_{n+1}$ matches the training data. Then, we implement stochastic gradient descent optimization to determine the parameters $\alpha^{n+1}$ and $\lambda^{n+1}$. Specifically, we define the loss function to be minimized as

$$
\begin{aligned}
F_{\alpha, \lambda}^{n+1} & :=\mathbb{E}\left[\left(p_{n+1}\left(S_{t_{n+1}}\right)-\tilde{p}\left(S_{t_{n+1}} \mid \mathcal{F}_{t_{n+1}}^M\right)\right)^{2}\right] \\
& =\mathbb{E}\left[\left(\sum_{k=1}^{K} \alpha_{k}^{n+1} \phi_{k}^{n+1}\left(S_{t_{n+1}}\right)-\tilde{p}\left(S_{t_{n+1}} \mid \mathcal{F}_{t_{n+1}}^M\right)\right)^{2}\right]
\end{aligned}
$$

Since the state $S_{t_{n+1}}$ is represented by spatial samples $\left\{\tilde{x}_{i}^{n+1}\right\}_{i=1}^{N}$, the fully calculated Monte Carlo simulation for the above loss function is given as

$$
F_{\alpha, \lambda}^{n+1} \approx \frac{1}{N} \sum_{i=1}^{N}\left(p_{n+1}\left(S_{t_{n+1}}=\tilde{x}_{i}^{n+1}\right)-\tilde{p}\left(S_{t_{n+1}}=\tilde{x}_{i}^{n+1} \mid \mathcal{F}_{t_{n+1}}^M\right)\right)^{2}
$$

Then, for pre-chosen initial estimates $\alpha^{n+1}(0)$ and $\lambda^{n+1}(0)$, we have the following stochastic gradient descent iteration to search for the parameters $\alpha^{n+1}$ and $\lambda^{n+1}$ :

\begin{subequations}
\label{Kernel:SGDupdate}
\begin{align}
\alpha^{n+1}(j+1)&=\alpha^{n+1}(j)-\left.\rho_{\alpha}^{j} \nabla_{\alpha} F_{\alpha, \lambda}^{n+1}\right|_{S_{t_{n+1}}=\bar{x}_{j}},  j=0,1,2, \cdots, J-1, \label{Kernel:SGDalpha} \\
\lambda^{n+1}(j+1)&=\lambda^{n+1}(j)-\left.\rho_{\lambda}^{j} \nabla_{\lambda} F_{\alpha, \lambda}^{n+1}\right|_{S_{t_{n+1}}=\bar{x}_{j}},  j=0,1,2, \cdots, J-1, \label{Kernel:SGDlambda}
\end{align}
\end{subequations}
where $\rho_{\alpha}^{j}$ and $\rho_{\lambda}^{j}$ are learning rates for the parameters $\alpha$ and $\lambda$, respectively, and $J$ is the total number of iterations corresponding to a stopping criteria. The gradients $\left.\nabla_{\alpha} F_{\alpha, \lambda}^{n+1}\right|_{S_{t_{n+1}}=\bar{x}_{j}}$ and $\left.\nabla_{\lambda} F_{\alpha, \lambda}^{n+1}\right|_{S_{t_{n+1}}=\bar{x}_{j}}$ of the cost function $F_{\alpha, \lambda}^{n+1}$ are single-sample representations of the gradients $\nabla_{\alpha} F_{\alpha, \lambda}^{n+1}$ and $\nabla_{\lambda} F_{\alpha, \lambda}^{n+1}$ by choosing a specific state sample $\bar{x}_{j}$ for the state variable $S_{t_{n+1}}$, where the sample $\bar{x}_{j}$ is picked among the sample set $\left\{\tilde{x}_{i}^{n+1}\right\}_{i=1}^{N}$. As a result, we improve the estimates for $\alpha^{n+1}$ and $\lambda^{n+1}$ gradually by comparing the kernel learned filtering density $p_{n+1}$ with
approximated filtering density values on samples $\left\{\tilde{x}_{i}^{n+1}\right\}_{i=1}^{N}$.

To use the approximated filtering density values more effectively, instead of picking samples uniformly from the sample set, we use importance sampling to choose samples according to their density values. In this way, it is more likely to consider higher filtering density values in the optimization procedure, which makes the stochastic gradient descent procedure more efficient. In this work, we let the covariance matrices for Gaussian type kernels be diagonal to reduce the dimension of optimization, and note that the scattered kernel centers also provide covariant features of the target filtering density.

\vspace{0.2cm}
\hspace{-0.72cm} \textbf{Resampling random spatial points}\vspace{0.2cm}

In the nonlinear filtering problem, the state equation is a diffusion process. Thus the random spatial samples propagated through the scheme (\ref{BSDE:ApproxCondexp6}) diffuse after several estimation steps. As a result, fewer and fewer spatial samples will remain in high probability regions of the filtering density as we estimate the target state step-by-step. So, we need to resample spatial points to make them better represent the filtering density by using the kernel learned updated filtering density $p_{n+1}$ to generate a set of new spatial samples, denoted by $\left\{x_{i}^{n+1}\right\}_{i=1}^{N}$, to replace the samples $\left\{\tilde{x}_{i}^{n+1}\right\}_{i=1}^{N}$, which follow the predicted filtering density. We want to point out that the kernel learned filtering density provides the conditional PDF for the target state over the entire state space. Therefore, this resampling procedure also allows us to consider probabilistically insignificant regions. On the other hand, the filtering density $p_{n+1}$ is essentially a combination of Gaussian kernels. Hence drawing samples from $p_{n+1}$ is very efficient. For example, to generate the sample $x_{i}^{n+1}$, we first use importance sampling to pick a kernel $\phi_{k}^{n+1}$ based on weights of kernels. Since $\phi_{k}^{n+1}$ is a Gaussian kernel, we can simply draw the sample $x_{i}^{n+1}$ from the Gaussian distribution $N\left(\hat{x}_{k}^{n+1}, \lambda_{k}^{n+1}\right)$.

\textbf{\textit{Remark}} The state spatial samples used in the backward SDE filter have similar behavior to the particles in the particle filter method, which roughly characterizes the filtering density in the nonlinear filtering problem. We want to emphasize that the backward SDE filter can also approximate the filtering density values on spatial samples. On the other hand, the particle filter only utilizes particle positions to construct an empirical distribution for the target state. Therefore, each spatial sample in the backward SDE filter carries more information about the state distribution than a particle in the particle filter. Moreover, since the kernel learned filtering density is a global continuous approximation for the state distribution, it covers a wide range in the state space, which can provide more robust/stable performance for the backward SDE filter.

\subsection{Stochastic optimization for control process}

In this subsection, we combine the numerical schemes for the adjoint FBSDEs system (\ref{PFSGD:EQ11})  and the BSDE filter algorithm to formulate an efficient stochastic optimization algorithm to solve for the optimal control process $u^{\ast}$.

On a time instant $t_n \in \Pi_{N_T}$, we have 
\begin{equation}
\label{PFSGD:EQ--13}
E[(J^{\ast})_u^{'}(u_{t_i}^{l,M}|_{t_{n}})|\mathcal{F}_{t_{n}}^M ] = \int_{\mathbb{R}^d} E[\Psi({t_i},S_{t_i},u_{t_i}^{l,M}|_{t_{n}})|S_{t_{n}}=x] \cdot p(x|\mathcal{F}_{t_{n}}^M) dx,
\end{equation}

where $t_i \geq t_n$ is a time instant after $t_n$.

Then, we use the approximate solutions $(\Bar{Y}_i,\Bar{Z}_i)$ of FBSDEs from schemes (\ref{PF-FBSDE:EQ15}) to replace  $(\Bar{Y}_{t_i},\Bar{Z}_{t_i})$ and the conditional distribution $p(S_{t_{n}}|\mathcal{F}_{t_{n}}^M)$ is approximated by the filtering density $p (S_{t_{n}}|\mathcal{F}_{t_{n}}^M)$ obtained from the BSDE filter through the kernel learning method. Then, we can generate S samples from $p (S_{t_{n}}|\mathcal{F}_{t_{n}}^M)$ and solve the optimal control $u_{t_n}^{\ast}$ through the following gradient descent optimization iteration
\begin{equation}
\label{PFSGD:EQ--18}
\begin{aligned}
u_{t_i}^{l+1,M}|_{t_n} =u_{t_i}^{l,M}|_{t_n} &-r\frac{1}{S}\sum_{s=1}^{S} E \Big[b_u^{\top}(t_i,S_{t_{i}},u_{t_i}^{l,M}|_{t_n}) \Bar{Y}_{{i}} +  \sigma_u^{\top}(t_i,S_{t_{i}},u_{t_i}^{l,M}|_{t_n})\Bar{Z}_{{i}} \\ 
&+f_u^{\top}(t_i,S_{t_{i}},u_{t_i}^{l,M}|_{t_n})| S_{t_n} = x_n^{(s)}  \Big]
\end{aligned}
\end{equation}
Then, the standard Monte Carlo method can approximate expectation $E \Big[\cdot | S_{t_n} = x_n^{(s)}  \Big]$ by $\Lambda$ samples:
\begin{equation}
\label{PFSGD:EQ--19}
\begin{aligned}
u_{t_i}^{l+1,M}|_{t_n} \approx u_{t_i}^{l,M}|_{t_n} &-r\frac{1}{S}\frac{1}{\Lambda}\sum_{s=1}^{S} \sum_{\lambda=1}^{\Lambda} \Big[b_u^{\top}(t_i,S_{t_{i}}^{(\lambda,s)},u_{t_i}^{l,M}|_{t_n}) \Bar{Y}_{{i}} +  \sigma_u^{\top}(t_i,S_{t_{i}}^{(\lambda,s)},u_{t_i}^{l,M}|_{t_n})\Bar{Z}_{{i}} \\ 
&+f_u^{\top}(t_i,S_{t_{i}}^{(\lambda,s)},u_{t_i}^{l,M}|_{t_n}) | S_{t_n} = x_n^{(s)} \Big]
\end{aligned}
\end{equation}

We can see from the above Monte Carlo approximation that in order to approximate the expectation in one gradient descent iteration step, we need to generate $S \times \Lambda$ samples. This is even more
computationally expensive when the controlled system is a high-dimensional process. 

So, we can also apply the idea of stochastic gradient descent (SGD) to improve the efficiency of classic gradient descent optimization as we used in the kernel learning for the BSDE filter. Instead of using the fully calculated Monte Carlo simulation to approximate the
conditional expectation, we use only one realization of $S_{t_i}$ to represent the expectation. For the conditional distribution of the controlled process, we can use the spatial samples to describe. So, we have
\begin{equation}
\label{PFSGD:EQ--0018}
\begin{aligned}
E[(J^{\ast})_u^{'}(u_{t_i}^{l,M}|_{t_{n}})|\mathcal{F}_{t_{n}}^M] \approx
b_u^{\top}(t_i,S_{t_{i}}^{(\hat{l},\hat{s})},u_{t_i}^{l,M}|_{t_n}) \Bar{Y}_{{i}}^{(\hat{l},\hat{s})} +  \sigma_u^{\top}(t_i,S_{t_{i}}^{(\hat{l},\hat{s})},u_{t_i}^{l,M}|_{t_n})\Bar{Z}_{{i}}^{(\hat{l},\hat{s})}
+f_u^{\top}(t_i,S_{t_{i}}^{(\hat{l},\hat{s})},u_{t_i}^{l,M}|_{t_n})
\end{aligned}
\end{equation}
where $l$ is the iteration index, the index $\hat{l}$ indicates that the random generation of the controlled process varies among the gradient
descent iteration steps .$S_{t_{i}}^{(\hat{l},\hat{s})}$ indicates a randomly generated realization of the controlled process with a randomly selected initial state $S_{t_{n}}^{(\hat{l},\hat{s})} = x_n^{\hat{s}}$ from the spatial sample cloud $\{x_n^{(s)}\}_{s=1}^S$.

Then, we have the following SGD schemes:
\begin{equation}
\label{PFSGD:EQ18}
\begin{aligned}
u_{t_i}^{l+1,M}|_{t_n} =u_{t_i}^{l,M}|_{t_n} &-r \Big(b_u^{\top}(t_i,S_{t_{i}}^{(\hat{l},\hat{s})},u_{t_i}^{l,M}|_{t_n}) \Bar{Y}_{{i}}^{(\hat{l},\hat{s})} +  \sigma_u^{\top}(t_i,S_{t_{i}}^{(\hat{l},\hat{s})},u_{t_i}^{l,M}|_{t_n})\Bar{Z}_{{i}}^{(\hat{l},\hat{s})} \\ 
&+f_u^{\top}(t_i,S_{t_{i}}^{(\hat{l},\hat{s})},u_{t_i}^{l,M}|_{t_n}) \Big)
\end{aligned}
\end{equation}
 $\Bar{Y}_{{i}}^{(\hat{l},\hat{s})}$ is the approximate solution $\Bar{Y}_i$ corresponding to the random sample $S_{t_{i}}^{(\hat{l},\hat{s})}$. And the path of $S_{t_{i}}^{(\hat{l},\hat{s})}$ is generated as following

\begin{equation}
\label{PFSGD:EQ19}
S_{t_{i+1}}^{(\hat{l},\hat{s})} = S_{t_{i}}^{(\hat{l},\hat{s})} + b(t_{i},S_{t_{i}}^{(\hat{l},\hat{s})},u_{t_i}^{l,M}|_{t_n})\triangle t_i + \sigma (t_{i},S_{t_{i}}^{(\hat{l},\hat{s})},u_{t_i}^{l,M}|_{t_n}) \sqrt{\triangle t_i} \omega_i^{(\hat{l},\hat{s})}    
\end{equation}
where $\omega_i^{(\hat{l},\hat{s})} \sim N(0,1)$. Then, an estimate for our desired data driven optimal
control at time instant $t_n$ is
\begin{equation*}
\hat{u}_{t_n} :=   u_{t_n}^{L,M}|_{t_n}
\end{equation*}
The scheme for FBSDE is 
\begin{equation}
\label{PFSGD:EQ20}
\begin{aligned}
&\Bar{Y}_{{i}}^{(\hat{l},\hat{s})} = \Bar{Y}_{{i+1}}^{(\hat{l},\hat{s})} + \Big[b_x(t_{i+1},S_{t_{i+1}}^{(\hat{l},\hat{s})},u_{t_i}^{l,M}|_{t_n})^{\top}\Bar{Y}_{{i+1}}^{(\hat{l},\hat{s})} + \sigma_x(t_{i+1},S_{t_{i+1}}^{(\hat{l},\hat{s})},u_{t_i}^{l,M}|_{t_n})^{\top}\Bar{Z}_{{i+1}}^{(\hat{l},\hat{s})}\\ 
& \;\;\;\;\;\;\;\;\;\;\;\;  + f_x(t_{i+1},S_{t_{i+1}}^{(\hat{l},\hat{s})},u_{t_i}^{l,M}|_{t_n})^{\top}\Big] \triangle t_i\\
&\Bar{Z}_{{i}}^{(\hat{l},\hat{s})} =\frac{\Bar{Y}_{{i+1}}^{(\hat{l},\hat{s})}\omega_i^{(\hat{l},\hat{s})}  }{\sqrt{\triangle t_i}}
\end{aligned}
\end{equation}
Then, we summarize the algorithm in the table below: 

\begin{algorithm}[]
\caption{BSDE-SGD algorithm for data driven feedback control problem}
\label{alg2}
Initialize the spatial sample cloud $\left\{x_{i}^{0}\right\}_{i=1}^{N} \sim p_{0}$, the number of kernels $K$, the learning rates $\rho_{\alpha}$ and $\rho_{\lambda}$, the number of iterations $L, J \in \mathbb{N}$, the total number of time steps $N_{T}$.
\begin{algorithmic}
\WHILE{$n=0,1,2,\cdots,N_T,$}
    \STATE Initialize an estimated control process $\{u_{t_i}^{(0,M)}|_{t_n}\}_{i=n}^{N_T}$ and a step-size $\rho$
    \FOR{SGD iteration steps $l=0,1,2,\cdots,L$,}
        \STATE \emph{\small$\bullet$} Simulate one realization of controlled process $\{S_{t_{i+1}}^{(\hat{l},\hat{s})}|_{t_n}\}_{i=n}^{N_T-1}$ through scheme (\ref{PFSGD:EQ19}) with $\{S_{t_{n}}^{(\hat{l},\hat{s})}\}=x_n^{(\hat{s})} \in \{x_n^{(s)}\}_{s=1}^{S}$;
        \STATE \emph{\small$\bullet$} Calculate solution $\{\Bar{Y}_{t_{i}}^{(\hat{l},\hat{s})}|_{t_n}\}_{i=N_T}^{n}$ of the FBSDEs system (\ref{PFSGD:EQ11}) corresponding to $\{S_{t_{i+1}}^{(\hat{l},\hat{s})}|_{t_n}\}_{i=n}^{N_T-1}$ through schemes (\ref{PFSGD:EQ20});
        \STATE \emph{\small$\bullet$} Update the control process to obtain $\{u_{t_i}^{(l+1,M)}|_{t_n}\}_{i=n}^{N_T}$ through scheme (\ref{PFSGD:EQ18});
    \ENDFOR
    \STATE \textbf{Update PDF of the state}
    \STATE \emph{\small$\bullet$} Propagate samples $\left\{x_{i}^{n}\right\}_{i=1}^{N}$ through the scheme (\ref{BSDE:ApproxCondexp6}) to get $\left\{\tilde{x}_{i}^{n+1}\right\}_{i=1}^{N}$.
    \STATE \emph{\small$\bullet$} Let $Y_{n}^{M_{t_{n}}}=p_{n}$ be the previous kernel learned filtering density. Solve the timeinverse FBSDEs system (\ref{BSDE:FBSDE2}) for $Y_{n+1}^{M_{t_{n}}}$ on spatial samples $\left\{\tilde{x}_{i}^{n+1}\right\}_{i=1}^{N}$ through the iterative scheme (\ref{BSDE:ApproxCondexp3}). $Y_{n+1}^{M_{t_{n}}}$ is the approximation for the predicted filtering density $\tilde{p}\left(S_{t_{n+1}} \mid \mathcal{F}_{t_{n}}^M\right)$.
    \STATE \emph{\small$\bullet$} Incorporate the observational information through Bayesian inference to get the updated filtering density $\tilde{p}\left(S_{t_{n+1}} \mid \mathcal{F}_{t_{n+1}}^M\right)$ on the spatial samples $\left\{\tilde{x}_{i}^{n+1}\right\}_{i=1}^{N}$.
    \STATE \emph{\small$\bullet$} Select kernel centers from the spatial samples $\left\{\tilde{x}_{i}^{n+1}\right\}_{i=1}^{N}$ by using the updated filtering density values on those samples.
    \STATE \emph{\small$\bullet$} Consider the approximated filtering density values $\left\{\tilde{p}\left(S_{t_{n+1}}=\tilde{x}_{i}^{n+1} \mid \mathcal{F}_{t_{n+1}}^M\right)\right\}_{i=1}^{N}$ as training data. Use the optimization procedure (\ref{Kernel:SGDupdate}) to obtain the kernel learned filtering density $p_{n+1}$ introduced in (\ref{Kernel:density});
    \STATE \emph{\small$\bullet$} Carry out the resampling procedure to generate new samples $\left\{x_{i}^{n+1}\right\}_{i=1}^{N}$ that follow the kernel learned filtering density $p_{n+1}$;
\ENDWHILE
\end{algorithmic}
\end{algorithm}

\section{Numerical examples}
\subsection{10D Linear Quadratic Experiment}
In this experiment, we want to show that the BSDE-SGD method can be applied to solve the classic linear quadratic control problem.

Assume A, B, R, K, Q are symmetric, positive definite. The forward process $Y$ and the observation process $M$ is given by
\begin{equation}
    \label{eq:turestate} 
    \begin{aligned}
        dY(t)&= A(u(t)-r(t))dt +\sigma B u(t)dW_t\\
        dM(t) &= sin(Y(t))dt + dB_t
    \end{aligned}
\end{equation}
The cost functional is given by
\begin{equation}
    \label{eq:cost} 
    J[u]=\frac{1}{2} E \left[ \int_0^T \langle R(Y_t-Y^{\ast}_t),(Y_t-Y^{\ast}_t) \rangle dt    + \frac{1}{2}\int_0^T\langle Ku_t,u_t \rangle dt  + \frac{1}{2} \langle QY_T,Y_T \rangle \right]
\end{equation}
we want to find $J(u^{\ast}) = \inf_{u \in \mathcal{U}_{ad} [0,T]} J(u)$.

\subsubsection{Construction of exact solutions}
An interesting fact of such example is that one can construct a time deterministic exact solution which depend
only on $x_0$.

By simplifying (\ref{eq:cost}), we have
\begin{equation}
    \label{eq:J1} 
    J[u]=\frac{1}{2} \int_0^T( RE [Y_t^T Y_t] - 2RY^{\ast T}_t E[Y_t] + RY^{\ast T}Y^{\ast} +\langle Ku_t,u_t \rangle dt) + \frac{1}{2}E[\langle QY_T,Y_T \rangle]
\end{equation}
Then, we define:
\begin{equation}
\begin{aligned}
\label{eq:J2} 
    S_t &= E[Y_t]=E[Y_0 + \int_0^T A(u(s)-r(s)) ds + \int_0^T \sigma B u(s)dW_s]\\
    &=E[Y_0 + \int_0^T A(u(s)-r(s) ds]
\end{aligned}
\end{equation}

Hence, we see that
\begin{equation}
\begin{aligned}
\label{eq:J3} 
    E [Y_t^T Y_t] &= E[(Y_0 + \int_0^t A(u(s)-r(s)) ds + \int_0^t \sigma B u(s)dW_s)^2]\\
    &=E[Y_0^T Y_0 + \int_0^t (u(s)-r(s))^TA^TA(u(s)-r(s) ds + Y_0^T\int_0^t A(u(s)-r(s)) ds \\
    &+ \int_0^t (u(s)-r(s))^TA^T Y_0 ds ] + E[\int_0^t \sigma^2u(s)^TB^TBu(s)ds]\\
    &=S_t^TS_t + \sigma^2\int_0^t u(s)^TB^TBu(s)ds
\end{aligned}
\end{equation}
And (\ref{eq:J3}) is true because all the terms are deterministic in time given $x_0$. Also, we observe that
\begin{equation}
\begin{aligned}
\label{eq:J4} 
    E [Y_T^T Y_T] &= E[(Y_0 + \int_0^T A(u(s)-r(s)) ds + \int_0^T \sigma B u(s)dW_s)^2]\\
    &=S_T^TS_T + \sigma^2\int_0^T u(s)^TB^TBu(s)ds
\end{aligned}
\end{equation}

As a result, we see that now (\ref{eq:J1}) takes the form:
\begin{equation}
\begin{aligned}
\label{eq:J5} 
    J[u]=&\frac{1}{2} \int_0^T(S_s^TRS_s - 2S_s^TRS_s^{\ast} +S_s^{\ast T}RS_s^{\ast} + u_s^T(\sigma^2 B^TQB+K)u_s)ds \\&+\frac{1}{2} \sigma^2 \int_0^T \int_0^t u_s^TB^TRBu_s dsdt + \frac{1}{2} S_T^TQS_T
    \end{aligned}
\end{equation}
by doing a simple integration by part, we have
\begin{equation}
\begin{aligned}
\label{eq:J6} 
    J[u]=&\frac{1}{2} \int_0^T(S_s^TRS_s - 2S_s^TRS_s^{\ast} +S_s^{\ast T}RS_s^{\ast} + u_s^T(\sigma^2 B^TQB+K)u_s)ds \\&+\frac{1}{2} \sigma^2 \int_0^T (T-t) u_s^TB^TRBu_s ds + \frac{1}{2} S_T^TQS_T
    \end{aligned}
\end{equation}
As a result, we have the following standard deterministic control problem:
\begin{equation}
\begin{aligned}
\label{eq:J7} 
    J[u]=&\frac{1}{2} \int_0^T \underbrace{(S_s^TRS_s - 2S_s^TRS_s^{\ast} +S_s^{\ast T}RS_s^{\ast} + u_s^T(\sigma^2 B^TQB+K)u_s + \sigma^2 (T-t) u_s^TB^TRBu_s)}_{2f} ds \\&+ \frac{1}{2} S_T^TQS_T
    \end{aligned}
\end{equation}
\begin{equation}
\label{eq:J8} 
\frac{dX_t}{dt}=  \underbrace{A(u(t)-r(t))}_{b}, X_{t_0} = X_0
\end{equation}
Then, one can form the following Hamiltonian
\begin{equation}
\label{eq:J9} 
H(x,p,u)=  bp + f
\end{equation}
where $p$ is $\nabla_x v$, $v$ is the value function.

Then, to find the optimal control, we have
\begin{align}
\label{eq:J10}
\frac{\partial}{\partial u} H = 0
\end{align}
which is 
\begin{align}
\label{eq:J11}
Ap + u(\sigma^2B^TRB(T-t)+(K+\sigma^2B^TQB)) = 0
\end{align}
So, we got
\begin{align}
\label{eq:J12}
u_t = \frac{-Ap(t) }{\sigma^2B^TRB(T-t)+(K+\sigma^2B^TQB)}
\end{align}
Also, notice that
\begin{align}
\label{eq:J13}
\frac{d }{dt}p(t) = -R(X_t-X_t^{\ast}), p(T) = QX_T
\end{align}
then,
\begin{align}
\label{eq:J14}
p(t) = QX_T + \int_t^T R(X_s-X_s^{\ast})ds
\end{align}
Combine (\ref{eq:J8}),(\ref{eq:J12}),(\ref{eq:J14}) together, we can solve the control of the system.

Set 
\begin{align}
\label{eq:J15}
    \mathbf A = \begin{bmatrix} 1 & 0.2 & \cdots & 0.2 & 0.2 \\  0.2 & 1 & \cdots & 0.2 & 0.2 \\  \vdots & \vdots & \ddots & \vdots & \vdots \\ 0.2 & 0.2 & \cdots & 1 & 0.2 \\  0.2 & 0.2 & \cdots & 0.2 & 1 \end{bmatrix} _{10 \times 10}
\end{align}
Set B, R, K ,Q be identity matrices.
To solve (\ref{eq:J14}), let 
\begin{equation*}
\label{eq:J17}
    \begin{aligned}
        X_t^1-X_t^{\ast 1} &:= t               &\;\; &X_t^2-X_t^{\ast 2} := cos(t) \\
        X_t^3-X_t^{\ast 3} &:= t^{2}           &\;\; &X_t^4-X_t^{\ast 4} := -2\pi sin(2\pi t)  \\
        X_t^5-X_t^{\ast 5} &:= -\frac{1}{2} e^{(T-t)}   & \;\;     &X_t^6-X_t^{\ast 6}:= \frac{1}{1+t} \\
        X_t^7-X_t^{\ast 7} &:= sec^2 (t)       & \;\;&X_t^8-X_t^{\ast 8} := \frac{1}{1+t^2} \\ 
        X_t^9-X_t^{\ast 9} &:= 1               & \;\;&X_t^{10}-X_t^{\ast {10}} := e^{(t-T)}
    \end{aligned}
\end{equation*}
Then we have
\begin{align}
\label{eq:J18}
    p(t) = \begin{bmatrix} X_t^1 \\X_t^2 \\X_t^3 \\ X_t^4\\ X_t^5\\ X_t^6\\ X_t^7\\ X_t^8\\ X_t^9\\ X_t^{10}\end{bmatrix} - \begin{bmatrix} \frac{T^2}{2}-\frac{t^2}{2} \\ sin(T)-sin(t) \\  \frac{T^3}{3}-\frac{t^3}{3}\\  cos(2\pi T)-cos(2\pi t) \\sinh(T)-sinh(t) \\log(1+T)-log(1+t) \\tan(T)-tan(t) \\arctan(T)-arctan(t) \\T-t \\1-e^{(t-T)} \end{bmatrix}
\end{align}

Let $\hat{r}(t)$ as 
\begin{equation*}
\label{eq:J19}
    \begin{aligned}
        & \hat{r}_t^1= \frac{-t^2/2}{\beta_t},  
        \hat{r}_t^2= \frac{-sin(t)}{\beta_t},  
        \hat{r}_t^3= \frac{-t^3}{3\beta_t}, 
        \hat{r}_t^4= \frac{-cos(2\pi t)}{\beta_t}, 
        \hat{r}_t^5= \frac{-sinh(t)}{\beta_t} \\
        & \hat{r}_t^6= \frac{-log(1+t)}{\beta_t},  
        \hat{r}_t^7= \frac{-tan(t)}{\beta_t}, 
        \hat{r}_t^8= \frac{-arctan(t)}{\beta_t}, 
        \hat{r}_t^9= \frac{-t}{\beta_t}, 
        \hat{r}_t^{10}= \frac{-e^{(t-T)}}{\beta_t} 
    \end{aligned}
\end{equation*}
where $ \beta_t = (1+\sigma^2)+\sigma^2 (T-t) $.
then 
$$r(t)=A\hat{r}(t)$$

Then, plug (\ref{eq:J18}) into (\ref{eq:J12}), then we solve (\ref{eq:J8})
\begin{align}
\label{eq:J20}
\frac{dX_t}{dt} &=  A(u(t)-r(t))\\ 
&=A(-\frac{Ap(t)}{\beta_t}- A\hat{r}(t)))
\end{align}
we get
\begin{align}
\label{eq:J21}
X_t = \begin{bmatrix} X_t^1 \\X_t^2 \\X_t^3 \\ X_t^4\\ X_t^5\\ X_t^6\\ X_t^7\\ X_t^8\\ X_t^9\\ X_t^{10} \end{bmatrix} = \frac{\alpha_t A^2}{\sigma^2 }(\begin{bmatrix} \frac{T^2}{2} \\sin(T) \\\frac{T^3}{3} \\ cos(2\pi T) \\sinh(T) \\log(1+T) \\tan(T) \\arctan(T) \\T \\1  \end{bmatrix}- \begin{bmatrix} X_T^1 \\X_T^2 \\X_T^3 \\ X_T^4 \\ X_T^5\\ X_T^6\\ X_T^7\\ X_T^8\\ X_T^9\\ X_T^{10}\end{bmatrix})
\end{align}
So, replace $X_t$ with $X_t^{\ast}$, we get 

\begin{align}
X_t^{\ast} =  \begin{bmatrix} X_t^{1,\ast} \\X_t^{2,\ast} \\X_t^{3,\ast} \\ X_t^{4,\ast} \\X_t^{5,\ast} \\X_t^{6,\ast}  \\X_t^{7,\ast} \\X_t^{8,\ast}  \\X_t^{9,\ast} \\X_t^{10,\ast}  \end{bmatrix}  =\begin{bmatrix} t \\cos(t) \\t^2 \\ -2 \pi sin(2\pi t) \\-\frac{1}{2} e^{(T-t)} \\ \frac{1}{1+t} \\ sec^2 (t) \\ \frac{1}{1+t^2} \\ 1 \\ e^{(t-T)} \end{bmatrix} + \frac{\alpha_t A^2}{\sigma^2 }(\begin{bmatrix} \frac{T^2}{2} \\sin(T) \\\frac{T^3}{3} \\ cos(2\pi T) \\sinh(T) \\log(1+T) \\tan(T) \\arctan(T) \\T \\1 \end{bmatrix} - \begin{bmatrix} X_T^1 \\X_T^2 \\X_T^3 \\ X_T^4 \\ X_T^5\\ X_T^6\\ X_T^7\\ X_T^8\\ X_T^9\\ X_T^{10} \end{bmatrix})
\end{align}
where $X_T^i$ can be obtained from the system (\ref{eq:J21}) by letting $t=T$, and 
$$\alpha_t = ln\frac{1+\sigma^2+\sigma^2 T}{(1+\sigma^2)+\sigma^2(T-t)}$$

Then, to find an the exact form by following the trajectory of $y_t$ in this setup, one will have to solve the following Coupled forward-backward ODE.

\begin{align}
	\frac{d X_t}{dt }&=  A(u_t-r_t) \ \ \ \ \ x_{n}=y_{t_n}\\ 
	\frac{d p(t)}{dt }&=X_t-X^*_t, \ \ \ p_T=x^{t_n, y_{t_n}}_T
\end{align}
with $u_t =-A p_t \big/\big( \sigma^2(T-t)+(1+\sigma^2) \big)$. 
As a result, we have
\begin{align}
	\frac{d X_t}{dt }&=  A(-Ap_t \big/\big( \sigma^2(T-t)+(1+\sigma^2) \big)-r_t) \ \ \ \ \ x_{n}=y_{t_n}\\ 
	\frac{d p(t)}{dt }&=X_t-X^*_t, \ \ \ p_T=x^{t_n, y_{t_n}}_T
\end{align}
That is, we need to solve the above coupled FBODE. Then, seeing that $p_t=x^{t_n, y_{t_n}}_T + \int^T_{t_n} X_s-X^*_s ds$. Writing $a_t:=1/\big( \sigma^2(T-t)+(1+\sigma^2) \big)$, we have 
\begin{align}\label{fbode}
	\frac{d X_t}{dt }&=A\big(-a_t A X_T-a_t A\int^T_{t_n} (X_s-X^*_s) ds -r_t \big),  \ \ \ \ \ x_{t_n}=y_{t_n}
\end{align}
To solve \eqref{fbode} numerically, we do a numerical discretization: 
\begin{align}
	x_{t_{n+1}}-x_{t_n}&=-a_{t_n}A^2X_T \Delta t-a_{t_n} (\Delta t)^2 A^2 \sum^{N-1}_{i=n}(X_{t_i}-X^*_{t_i}) - A r_t \nonumber \\ 
	 \Rightarrow  A r_t -a_{t_n} (\Delta t)^2 A^2\sum^{N-1}_{i=n} X^*_{t_i}&=X_{t_n}-X_{t_{n+1}}-a_{t_n} (\Delta t)^2 A^2 \sum^{N-1}_{i=n}X_{t_i}-a_{t_n}A^2X_T, \ \ \ x_{t_n}=y_{t_n} \label{fbode-num}
\end{align}
We can put \eqref{fbode-num} into a large linear system, and solve it numerically. 

\subsubsection{Result}
\begin{figure}[H]
    \centering
    \includegraphics[width=15cm]{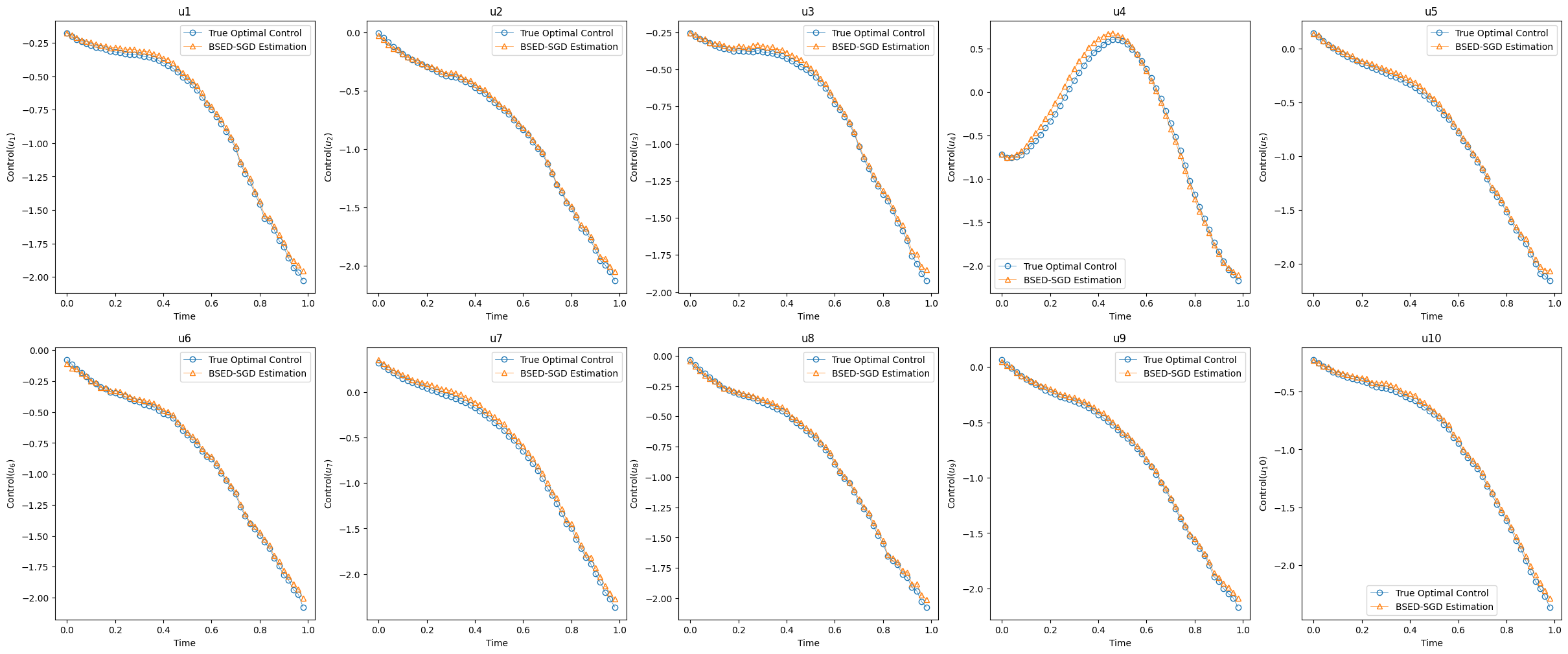} %\hfill
        \caption{Estimated control vs True optimal control}
        \label{fig1}
\end{figure}

\begin{figure}[H]
    \centering
    \includegraphics[width=15cm]{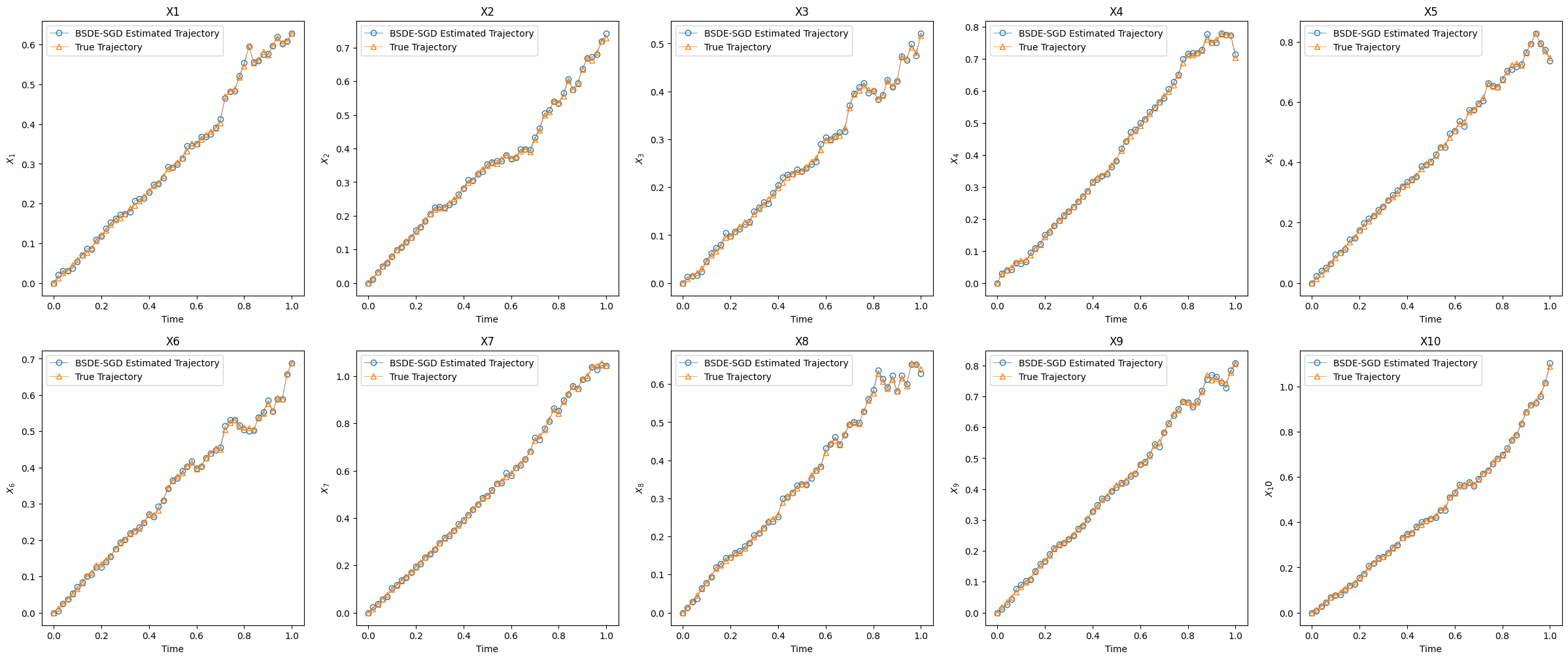} %\hfill
        \caption{Estimated trajectory vs True trajectory}
        \label{fig12}
\end{figure}

Set iteration $L=10^4, \sigma=0.1$, $T=1$, $X_0 = 0$ we have the following results in figure \ref{fig1} and \ref{fig12}. We demonstrate the effectiveness of our data driven feedback control algorithm
by comparing with the analytical solution, where the analytical optimal control is calculated
from (\ref{eq:J12}) with the exact controlled state and we let the analytical optimal control drive its
corresponding controlled system through (\ref{eq:turestate}) to get the analytical controlled process. The estimated optimal control and its corresponding controlled process are calculated by our data driven feedback control Algorithm based on fully nonlinear observations. figure \ref{fig1} shows that the estimated control is very close to the analytical optimal control. figure \ref{fig12} shows the controlled process governed by our data driven feedback control well aligns with the real state of the analytical controlled process.
\subsubsection{Testing convergence}
\begin{figure}[H]
\centering
\begin{subfigure}{.5\textwidth}
    \centering
    \includegraphics[width=.85\linewidth]{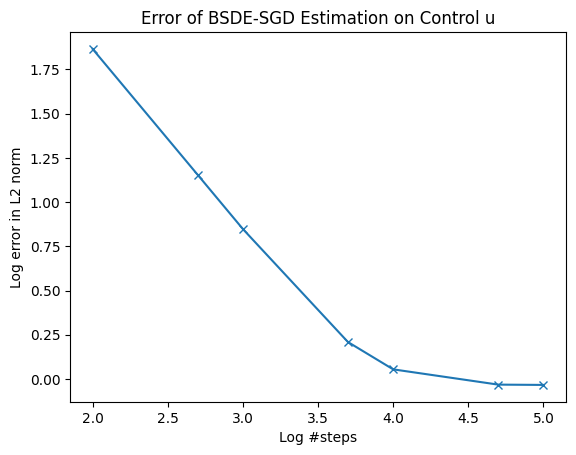} 
    \caption{Control u}
    \label{fig2}
\end{subfigure}%
\begin{subfigure}{.5\textwidth}
\centering
    \includegraphics[width=.85\linewidth]{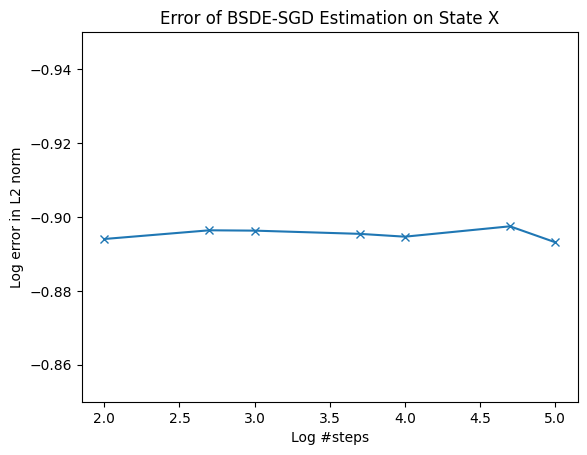} 
    \caption{State X }
    \label{fig3}
\end{subfigure}
\caption{Error of BSDE-SGD Estimation on Control u and State X }
\end{figure}

To test the convergence of our algorithm, we increase the number of iterations, we got the result in (figure \ref{fig2}), where the error between the estimated control and the analytical control is the average of $||u^{est}-u^{\ast}||_2$ of 50 repeated experiments. The result shows that when the number of iterations in the stochastic gradient descent algorithm increases, the error decreases and converges. This means we can get a better result by increasing the number of iterations. At the same time, we also got the result in (figure \ref{fig3}), where the error between the controlled process governed by our data driven feedback control and the analytical controlled process is the average of $||X^{est}-X^{ture \; trajectory}||_2$ of the same 50 repeated experiments. It tells that the BSDE filter provides a very accurate and stable estimation of the true state X.

\subsection{PDE Experiment}
In this experiment, we want to show that the BSDE-SGD method can be applied to solve the PDE with control terms. The idea is we use a well-studied scheme to discretize the PDE where we get a high-dimensional ODE system. Then we can add a diffusion term to convert it into a high-dimensional SDE system. We also assume that we do not have the true state and we can only have observation data. Finally, we can solve this stochastic optimal feedback control problem.

First, the problem is to find the optimal control of the 1-D heat equation:
\begin{equation}
    \label{eq:1Dheat} 
    \begin{aligned}
        & P_t  = aP_{xx} + bP_x + cu(x,t) + df(x,t), \;\;\;\; (x,t)\in [0,1] \times[0,1] \\
        & P(x,0) = 10x(1-x)\frac{10}{(1+0)} e^{-\frac{x-0}{(1+0)^2}} \\
        &f(x,t) = x(1-x)t(1-t) \\
        &P(0,t)  = 0\\
        &P(1,t)  = 0
    \end{aligned}
\end{equation}
The cost function J is defined by
\begin{equation}
    \label{eq:1Dheatcost} 
    J[u]=\frac{1}{2} E \left[ \int_0^T \langle RP(t),P(t) \rangle dt    + \int_0^T\langle Qu(t),u(t) \rangle dt  + \langle KP(T),P(T) \rangle \right]
\end{equation}
The discretized system is given by:

\begin{equation}
    \label{eq:dis1} 
        \frac{d}{dt} P  = AP + Bu + Cf
\end{equation}
using forward-time for $P_t$, central-space for $P_{xx}$ and backward-space for $P_x$, the matrices $A$,$B$, and $C$ are given by

\begin{align}
\label{eq:dis2}
    \mathbf A = \frac{a}{h^2}\begin{bmatrix} -2 & 1 &   &   &   \\  1 & -2 & 1 &   &   \\    & \ddots & \ddots & \ddots &   \\   &  & 1 & -2 & 1 \\    &   &   & 1 & -2 \end{bmatrix} _{n \times n} + \frac{b}{h}\begin{bmatrix} 1 &   &   &   &   \\  -1 & 1 &   &   &   \\    & \ddots & \ddots &   &   \\   &  & -1 & 1 &   \\    &   &   & -1 & 1 \end{bmatrix} _{n \times n}
\end{align}

\begin{align}
\label{eq:dis2_2}
    \mathbf B = \begin{bmatrix} 0 \\c_1  \\ \vdots  \\ c_{n-1} \\ 0\end{bmatrix} _{n \times 1}, \mathbf C = \begin{bmatrix} 0 \\d_1  \\ \vdots  \\ d_{n-1} \\ 0\end{bmatrix} _{n \times 1}
\end{align}

For the system of (\ref{eq:dis1}), the Riccati equations are given by:
\begin{equation}
    \label{eq:dis3} 
        \frac{d}{dt} G(t)  = -A^T G(t) - G(t)A - G(t)BR^{-1}B^TG(t) + Q
\end{equation}
The optimal control $u$ for the above stochastic optimal control problem can be derived
analytically as
\begin{equation}
    \label{eq:dis4} 
        u^{\ast}  = (-R^{-1}B^TG(t)P(t) - df(t))/c 
\end{equation}

To ensure the discretization is stable, we set $a=b= 5 \times 10^{-5}$. We also set $c = d =1$, $R$, $Q$, and $K$ be identity matrix.   

To apply the BSDE-SGD method, we can add the diffusion term to the system (\ref{eq:dis1}) and the observation operator is fully nonlinear:

\begin{equation}
    \label{eq:dis5} 
    \begin{aligned}
        dP_t  &= (AP + Bu + Cf) dt + \sigma dW_t \\
        dM_t &= sin(P_t)dt + \eta dB_t
    \end{aligned}
\end{equation}

Then we can use the BSDE-SGD method to the stochastic optimal feedback control problem giving the cost function J in (\ref{eq:1Dheatcost}).

\begin{figure}[H]
\centering
\begin{subfigure}{.33\textwidth}
    \centering
    \includegraphics[width=.9\linewidth]{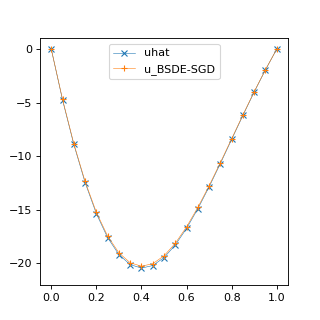} 
    \caption{t=0}
\end{subfigure}%
\begin{subfigure}{.33\textwidth}
    \centering
    \includegraphics[width=.9\linewidth]{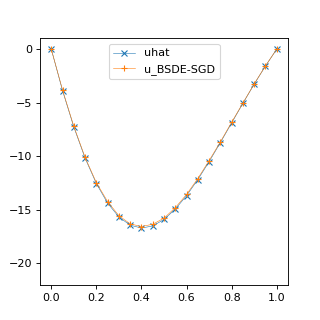} 
    \caption{t=0.2}
\end{subfigure}%
\begin{subfigure}{.33\textwidth}
    \centering
    \includegraphics[width=.9\linewidth]{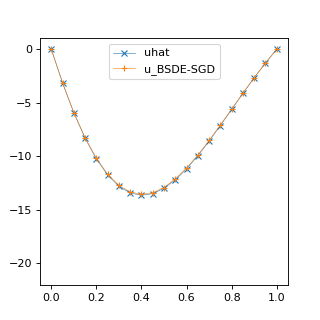} 
    \caption{t=0.4}
\end{subfigure}%

\begin{subfigure}{.33\textwidth}
    \centering
    \includegraphics[width=.9\linewidth]{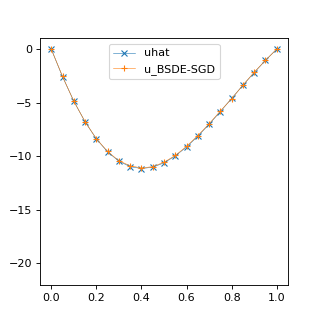} 
    \caption{t=0.6}
\end{subfigure}%
\begin{subfigure}{.33\textwidth}
    \centering
    \includegraphics[width=.9\linewidth]{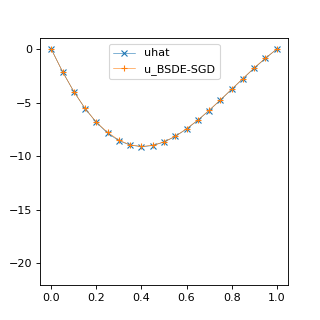} 
    \caption{t=0.8}
\end{subfigure}%
\begin{subfigure}{.33\textwidth}
    \centering
    \includegraphics[width=.9\linewidth]{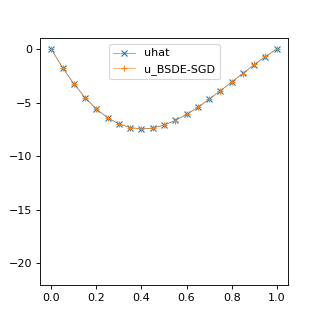} 
    \caption{t=1}
\end{subfigure}%
\caption{BSDE-SGD Estimation VS Optimal Control}
\label{HeatBSDEu}
\end{figure}

\begin{figure}[H]
\centering
\begin{subfigure}{.33\textwidth}
    \centering
    \includegraphics[width=.9\linewidth]{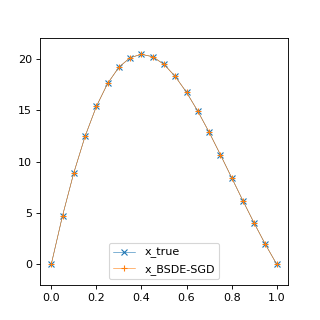} 
    \caption{t=0}
\end{subfigure}%
\begin{subfigure}{.33\textwidth}
    \centering
    \includegraphics[width=.9\linewidth]{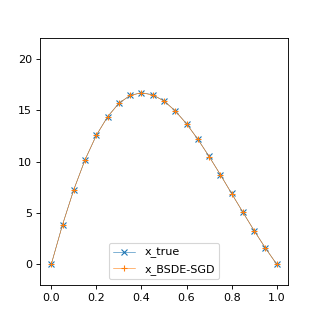} 
    \caption{t=0.2}
\end{subfigure}%
\begin{subfigure}{.33\textwidth}
    \centering
    \includegraphics[width=.9\linewidth]{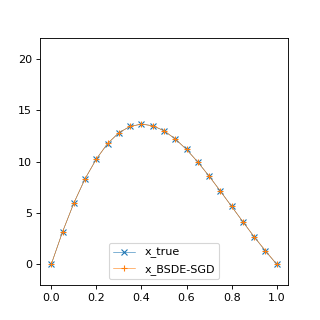} 
    \caption{t=0.4}
\end{subfigure}%

\begin{subfigure}{.33\textwidth}
    \centering
    \includegraphics[width=.9\linewidth]{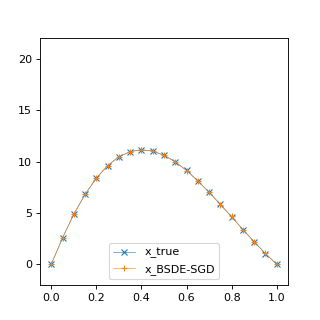} 
    \caption{t=0.6}
\end{subfigure}%
\begin{subfigure}{.33\textwidth}
    \centering
    \includegraphics[width=.9\linewidth]{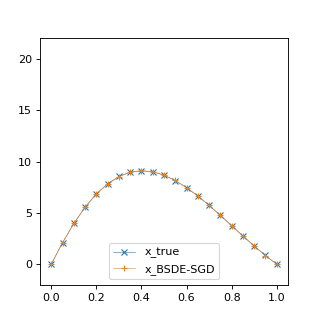} 
    \caption{t=0.8}
\end{subfigure}%
\begin{subfigure}{.33\textwidth}
    \centering
    \includegraphics[width=.9\linewidth]{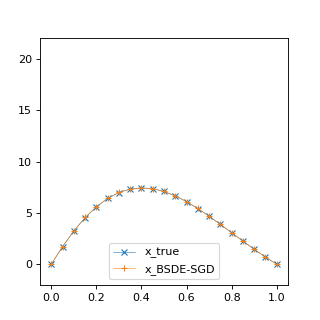} 
    \caption{t=1}
\end{subfigure}%
\caption{BSDE-SGD Estimation VS True State}
\label{HeatBSDEX}
\end{figure}

Set iteration $L=10^4, \sigma=0.1, \eta = 0.1$, and number of nodes on space $n=21$. We have the following results. From Figure \ref{HeatBSDEu}, we see that the estimated control is very accurate compared with the analytical solution of control. From Figure \ref{HeatBSDEX}, we also see that the estimated heat distribution governed by our data driven feedback control well aligns with the analytical real heat distribution. We demonstrate the effectiveness of our data driven feedback control algorithm of solving a PDE control problem by converting it into a high-dimensional SDE system governed by the control.

\subsection{3D Dubins Airplane maneuvering problem}

In this example, we solve a Dubins airplane maneuvering problem. The controlled
process is described by the following nonlinear controlled dynamics:
\begin{align}
	dS_t &=\begin{bmatrix} dX_t \\dY_t \\dZ_t\end{bmatrix} = \begin{bmatrix} v\cos(\theta_t)\cos(\phi_t) dt  + \sigma dW_{1t} \\ v\cos(\theta_t)\sin(\phi_t)dt   + \sigma dW_{2t} \\ v\sin(\theta_t)  + \sigma dW_{3t}  \end{bmatrix}\\  d \theta_t &= u_t + \sigma^2 dB_{1t} \\ d \phi_t &= p_t + \sigma^2 dB_{2t} \\
 dM_t &= [arctan(\frac{X_t+3}{Y_t+2}),arctan(\frac{X_t-2}{Y_t+2}),arctan(\frac{Z_t-2}{Y_t+2})]^T + \eta_t
	  \label{dubinsair}
\end{align}
where the pair $(X, Y, Z)$ gives the position of a airplane-like robot moving in the 3D plane. $\theta$ is
the angle that controls the moving direction of the airplane vertically, which is governed by the
control action $u_t$. $\phi$ is the angle that controls the moving direction of the airplane horizontally which is governed by the
control action $p_t$. $\sigma$ is the noise that perturbs the motion and control actions. Assume that we do not have direct observations on the robot. Instead, we use two detectors
located on different observation platforms at $(-3, -2, 2)$ and $(2, -2,-2)$ to collect the bearing angles of
the target robot as indirect observations. So, we have the observation process $M_t$ Given the expected path $S^{\ast}$, the airplane should follow it and arrive at the terminal position on time. The performance cost functional based on observational
data that we aim to minimize is defined as:
\begin{equation}
    \label{eq:cost2d} 
    J[u]= E \left[ \frac{1}{2}\int_0^T \langle R(S_t-S^{\ast}_t),(S_t-S^{\ast}_t) \rangle dt    + \frac{1}{2}\int_0^T\langle Ku_t,u_t \rangle dt  + \langle Q(S_T-S^{\ast}_T),(S_T-S^{\ast}_T) \rangle \right]
\end{equation}

In our numerical experiments, we let the airplane start from $(X_0,Y_0,Z_0)=(0,0.5,0 )$ to $(X_T,Y_T,Z_T)=(0,0.5,1)$. The expected path $S^{\ast}_t$ is 
$$X(t) = 0.5 \sin(2\pi t), \;
Y(t) = 0.5 \cos(2\pi t), \;
Z(t) = t$$
Other settings are $T=1$, $\Delta t  =0.02$ i.e. $N_T=50$, $\sigma=0.1$, $\eta_t \sim N(0,0.1)$, $L=1000$, $K=1$ and the initial heading direction is $\theta = 0 $ and $\phi = arctan(1 / 2\pi)$. To emphasize the importance of following the expected path and arriving at the target location at the terminal time, let $Q=20$ and $R=40$.
\begin{figure}[H]
    \centering
    \includegraphics[width=9cm]{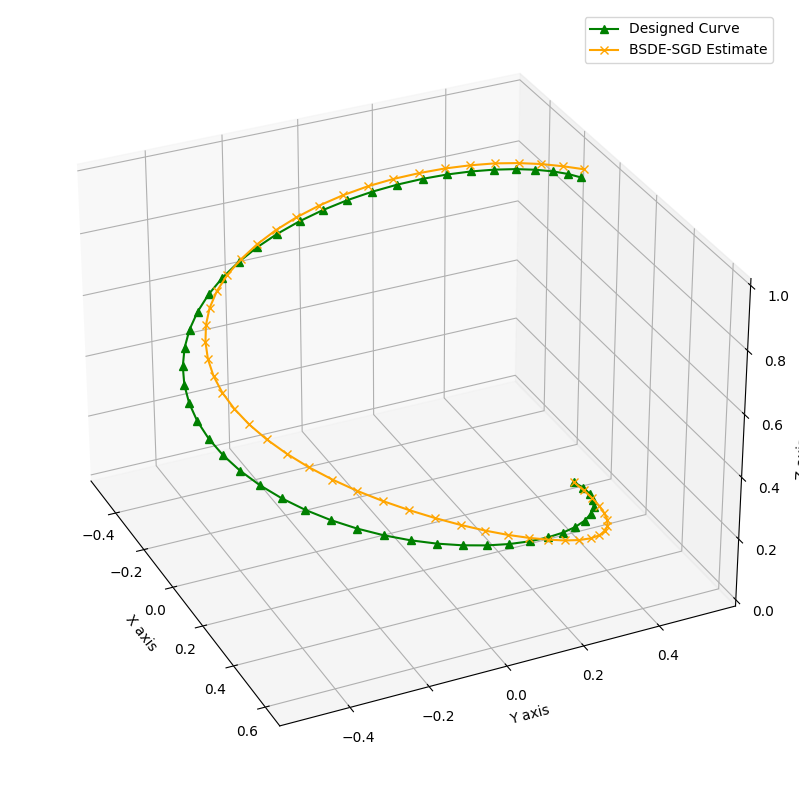} %\hfill
        \caption{BSDE-SGD Estimate VS Designed trajectory}
        \label{3Dair1}
\end{figure}

\begin{figure}[H]
    \centering
    \includegraphics[width=.85\linewidth]{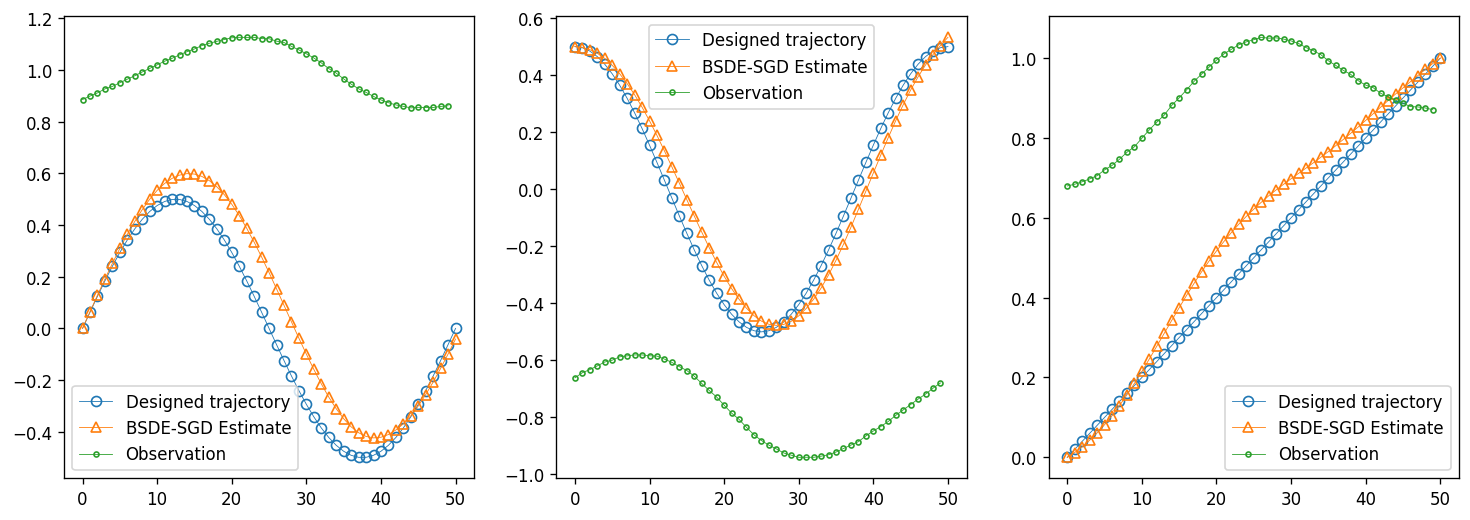} %\hfill
        \caption{BSDE-SGD Estimate VS Designed trajectory in Each Dimension (X, Y, Z from left to right)}
        \label{3Dair2}
\end{figure}
We can see that the observation is highly non-linear. But the estimation of the BSDE-SDG method gives the result where the airplane moves to the designed terminal position. At the meantime, it also follows the designed trajectory closely.

\subsubsection{Compare with the Particle Filter Method}
\begin{figure}[H]
    \centering
    \includegraphics[width=.85\linewidth]{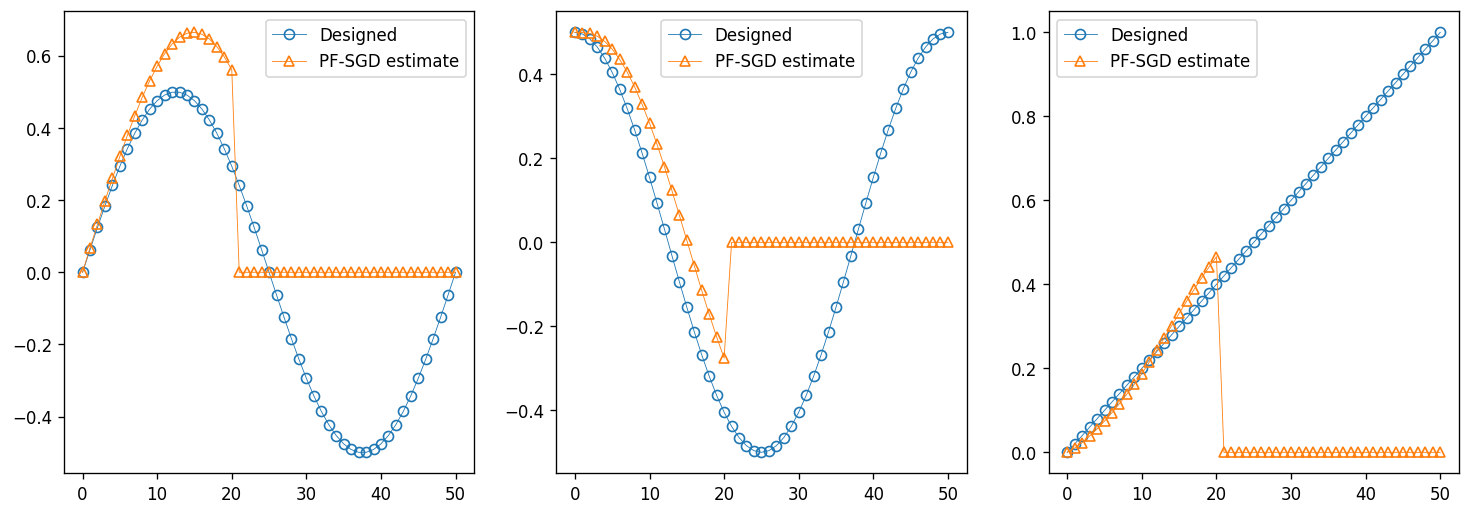} %\hfill
        \caption{PF-SGD Estimate VS Designed trajectory in Each Dimension (X, Y, Z from left to right)}
        \label{JUMP_PF_SDG}
\end{figure}

\begin{figure}[H]
    \centering
    \includegraphics[width=.85\linewidth]{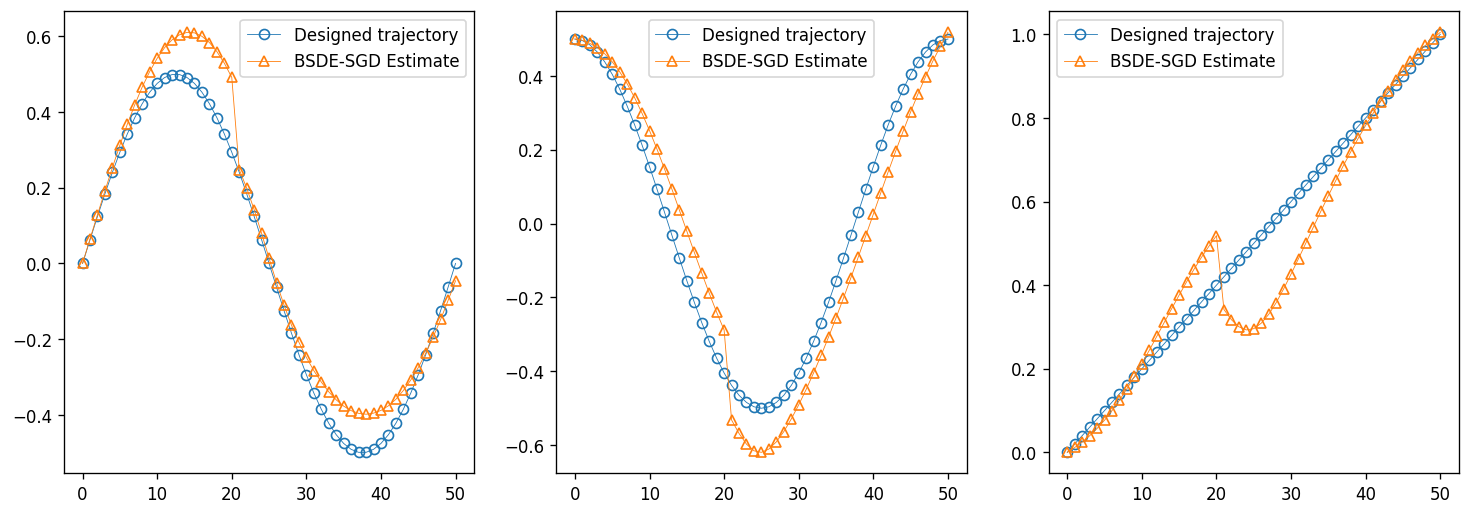} %\hfill
        \caption{BSDE-SGD Estimate VS Designed trajectory in Each Dimension (X, Y, Z from left to right)}
        \label{JUMP_BSDE_SDG}
\end{figure}

In this part, we will add a big jump and huge noise to the Airplane and the system to simulate the situation where the airplane is suddenly in front of a thundercloud. In this case, the airplane has to dive quickly to avoid crossing the cloud. In the meantime, radars on the ground will receive huge noise due to the weather. In the figure \ref{JUMP_PF_SDG}, we see after it happens, the particle filter method fails to track to position of the airplane. But figure \ref{JUMP_BSDE_SDG} shows that the BSDE filter works well to track the airplane when a big jump happens. It also shows that when the filter works well, our data driven feedback control can "drag" the airplane back to the designed path and move to the target position.

\subsubsection{Compare with the Dynamic Programming Method}
In this section, we will compare our algorithm with the Dynamic Programming method. In figure \ref{3Dair1}, we have shown the result of our algorithm where we used 1,000 iterations which cost 67 seconds. To compare with the Dynamic Programming method, we set the 3 by 3 grid for the control, $T=1$, $N_T=8$. Note that we already know the estimated control from our algorithm, so the 3 by 3 grid for the control is very close to its best choice. But even with this "trick", it still takes 627.5 seconds to get the result in figure \ref{DP} below.
\begin{figure}[H]
    \centering
    \includegraphics[width=.5\linewidth]{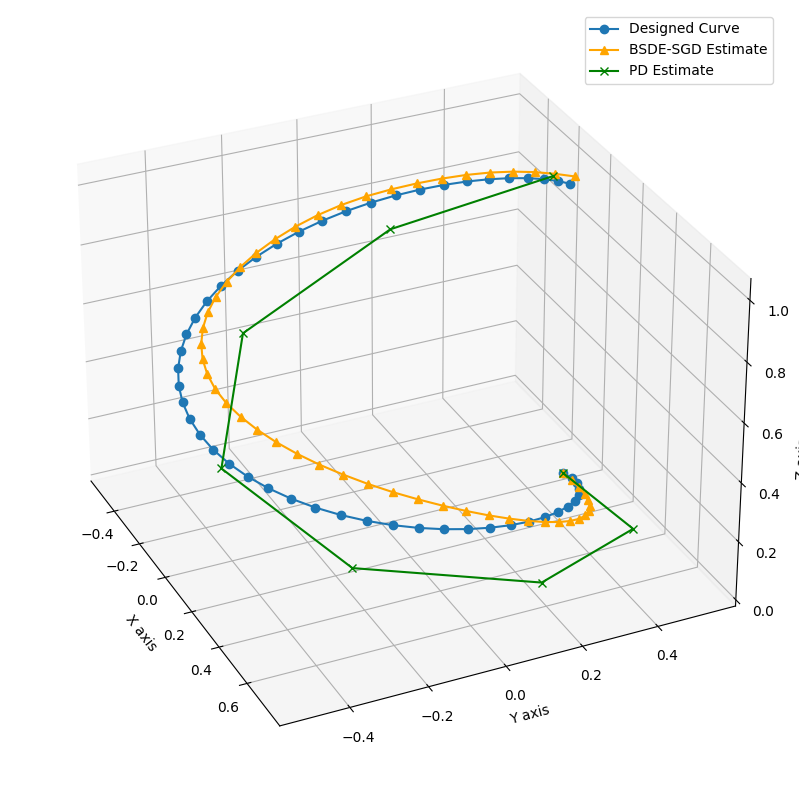} \hfill
\end{figure}

\begin{figure}[H]
    \centering
    \includegraphics[width=.85\linewidth]{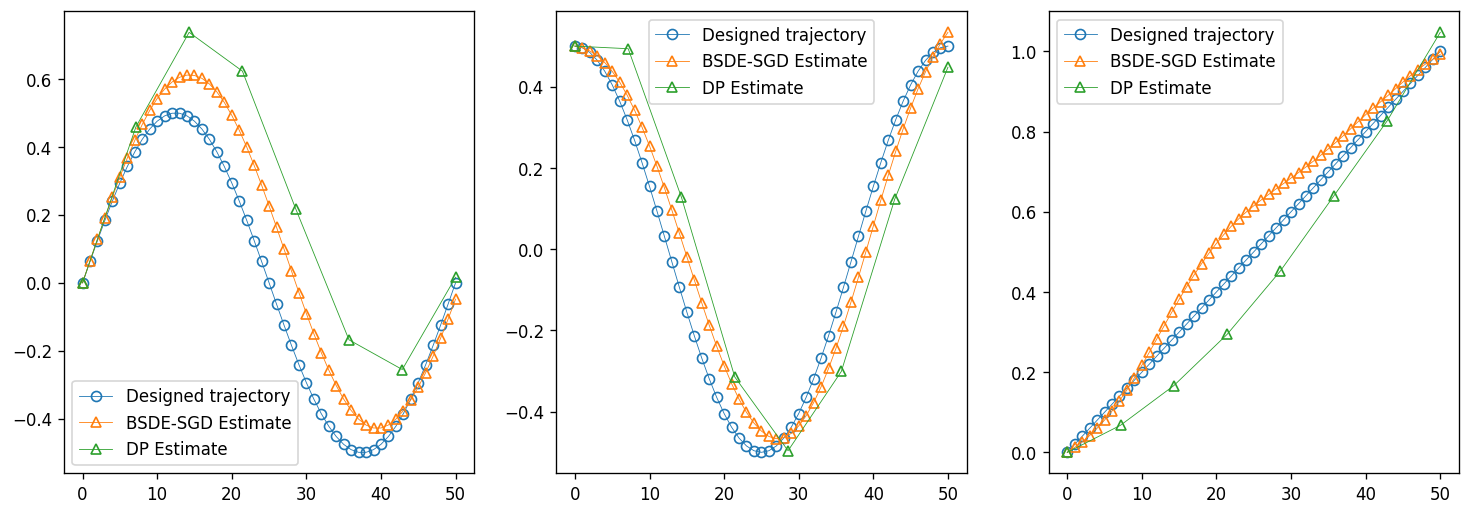} %\hfill
        \caption{DP Estimate VS Designed trajectory in Each Dimension (X, Y, Z from left to right)}
        \label{DP}
\end{figure}
From figure \ref{DP}, we see the Dynamic Programming method works, but the result is very rough. If we when to get a better result, we will need to increase the time partition and set a control grid that covers more points. But by just setting  $N_T=50$, it will take $6 \times 10^{42}$ seconds to finish.

%\newpage
%\renewcommand*{\bibname}{References}
%\printbibliography
\bibliographystyle{plain}
%\bibliography{myrefs, Reference, Reference_control}

\begin{thebibliography}{10}

\bibitem{MCMC-PF}
C.~Andrieu, A.~Doucet, and R.~Holenstein.
\newblock Particle markov chain monte carlo methods.
\newblock {\em J. R. Statist. Soc. B}, 72(3):269--342, 2010.

\bibitem{Bao_EAJAM20}
R.~Archibald, F.~Bao, and J.~Yong.
\newblock A stochastic gradient descent approach for stochastic optimal
  control.
\newblock {\em East Asian Journal on Applied Mathematics}, 10(4):635--658,
  2020.

\bibitem{Bao_Control_20}
R.~Archibald, F.~Bao, J.~Yong, and T.~Zhou.
\newblock An efficient numerical algorithm for solving data driven feedback
  control problems.
\newblock {\em Journal of Scientific Computing}, 85(51), 2020.

\bibitem{Bao_SINUM_SNN}
Richard Archibald, Feng Bao, Yanzhao Cao, and Hui Sun.
\newblock Numerical analysis for convergence of a sample-wise backpropagation
  method for training stochastic neural networks.
\newblock {\em SIAM Journal on Numerical Analysis}, 62(2):593--621, 2024.

\bibitem{Bao_RL}
Richard Archibald, Feng Bao, and Jiongmin Yong.
\newblock A stochastic maximum principle approach for reinforcement learning
  with parameterized environment.
\newblock {\em Journal of Computational Physics}, 488:112238, 2023.

\bibitem{archibald2020efficient}
Richard Archibald, Feng Bao, Jiongmin Yong, and Tao Zhou.
\newblock An efficient numerical algorithm for solving data driven feedback
  control problems.
\newblock {\em Journal of Scientific Computing}, 85(2):1--27, 2020.

\bibitem{Bao_CiCP20}
F.~Bao, Y.~Cao, and P.~Maksymovych.
\newblock Backward sde filter for jump diffusion processes and its applications
  in material sciences.
\newblock {\em Communications in Computational Physics}, 27:589--618, 2020.

\bibitem{Bao_AA20}
F.~Bao, Y.~Cao, and J.~Yong.
\newblock Data informed solution estimation for forward backward stochastic
  differential equations.
\newblock {\em Analysis and Applications, to appear}, 2020.

\bibitem{BaoC20142}
Feng Bao, Yanzhao Cao, and Xiaoying Han.
\newblock Forward backward doubly stochastic differential equations and optimal
  filtering of diffusion processes.
\newblock {\em Communications in Mathematical Sciences}, 18(3):635--661, 2020.

\bibitem{Bao_first}
Feng Bao, Yanzhao Cao, Amnon Meir, and Weidong Zhao.
\newblock A first order scheme for backward doubly stochastic differential
  equations.
\newblock {\em SIAM/ASA J. Uncertain. Quantif.}, 4(1):413--445, 2016.

\bibitem{bao2016first}
Feng Bao, Yanzhao Cao, Amnon Meir, and Weidong Zhao.
\newblock A first order scheme for backward doubly stochastic differential
  equations.
\newblock {\em SIAM/ASA Journal on Uncertainty Quantification}, 4(1):413--445,
  2016.

\bibitem{Bao_zakai}
Feng Bao, Yanzhao Cao, Clayton Webster, and Guannan Zhang.
\newblock A hybrid sparse-grid approach for nonlinear filtering problems based
  on adaptive-domain of the {Z}akai equation approximations.
\newblock {\em SIAM/ASA J. Uncertain. Quantif.}, 2(1):784--804, 2014.

\bibitem{bao2011numerical}
Feng Bao, Yanzhao Cao, and Weidong Zhao.
\newblock Numerical solutions for forward backward doubly stochastic
  differential equations and zakai equations.
\newblock {\em Visualization of Mechanical Processes: An International Online
  Journal}, 1(4), 2011.

\bibitem{BSDE_filter}
Feng Bao and Vasileios Maroulas.
\newblock Adaptive meshfree backward {SDE} filter.
\newblock {\em SIAM J. Sci. Comput.}, 39(6):A2664--A2683, 2017.

\bibitem{Bellman1957}
R.~Bellman.
\newblock {\em Dynamic Programming}.
\newblock 1957.

\bibitem{bottou2007tradeoffs}
L{\'e}on Bottou and Olivier Bousquet.
\newblock The tradeoffs of large scale learning.
\newblock {\em Advances in neural information processing systems}, 20, 2007.

\bibitem{Charalambous-98}
Charalambos~D. Charalambous and Robert~J. Elliott.
\newblock Classes of nonlinear partially observable stochastic optimal control
  problems with explicit optimal control laws.
\newblock {\em SIAM J. Control Optim.}, 36(2):542--578, 1998.

\bibitem{CT1}
A.~J. Chorin and X.~Tu.
\newblock Implicit sampling for particle filters.
\newblock {\em Proc. Nat. Acad. Sc. USA}, 106:17249--17254, 2009.

\bibitem{Do2}
A.~Doucet, N.~de~Freitas, and N.~Gordon.
\newblock {\em Sequential {M}onte {C}arlo Methods in Practice}.
\newblock Springer, New York, 2001.

\bibitem{Evense_EnKF}
G.~Evensen.
\newblock The ensemble {K}alman filter for combined state and parameter
  estimation: {M}onte {C}arlo techniques for data assimilation in large
  systems.
\newblock {\em IEEE Control Syst. Mag.}, 29(3):83--104, 2009.

\bibitem{Feng_HJB_2013}
Xiaobing Feng, Roland Glowinski, and Michael Neilan.
\newblock Recent developments in numerical methods for fully nonlinear second
  order partial differential equations.
\newblock {\em SIAM Rev.}, 55(2):205--267, 2013.

\bibitem{Fleming-Pardoux-1982}
Wendell~H. Fleming and \'{E}tienne Pardoux.
\newblock Optimal control for partially observed diffusions.
\newblock {\em SIAM J. Control Optim.}, 20(2):261--285, 1982.

\bibitem{GPM_2017}
Bo~Gong, Wenbin Liu, Tao Tang, Weidong Zhao, and Tao Zhou.
\newblock An efficient gradient projection method for stochastic optimal
  control problems.
\newblock {\em SIAM J. Numer. Anal.}, 55(6):2982--3005, 2017.

\bibitem{gordon1993novel}
Neil~J Gordon, David~J Salmond, and Adrian~FM Smith.
\newblock Novel approach to nonlinear/non-gaussian bayesian state estimation.
\newblock In {\em IEE proceedings F (radar and signal processing)}, volume 140,
  pages 107--113. IET, 1993.

\bibitem{particle-filter}
N.J Gordon, D.J Salmond, and A.F.M. Smith.
\newblock Novel approach to nonlinear/non-gaussian bayesian state estimation.
\newblock {\em IEE PROCEEDING-F}, 140(2):107--113, 1993.

\bibitem{Haussmann-1982}
U.~G. Haussmann.
\newblock On the existence of optimal controls for partially observed
  diffusions.
\newblock {\em SIAM J. Control Optim.}, 20(3):385--407, 1982.

\bibitem{hofmann2008kernel}
Thomas Hofmann, Bernhard Sch{\"o}lkopf, and Alexander~J Smola.
\newblock Kernel methods in machine learning.
\newblock 2008.

\bibitem{Kang-PF}
Kai Kang, Vasileios Maroulas, Ioannis Schizas, and Feng Bao.
\newblock Improved distributed particle filters for tracking in a wireless
  sensor network.
\newblock {\em Comput. Statist. Data Anal.}, 117:90--108, 2018.

\bibitem{Kloeden1992}
Peter~E. Kloeden and Eckhard Platen.
\newblock {\em Stochastic Differential Equations}, pages 103--160.
\newblock Springer Berlin Heidelberg, Berlin, Heidelberg, 1992.

\bibitem{Zhao_BSDE_Control_17}
T~Kong, W~Zhao, and T.~Zhou.
\newblock High order numerical schemes for second order fbsdes with
  applications to stochastic optimal control.
\newblock {\em Communications in Computational Physics}, 21:808--834, 2017.

\bibitem{X_Li_Drift}
Xin Li, Feng Bao, and Kyle Gallivan.
\newblock A drift homotopy implicit particle filter method for nonlinear
  filtering problems.
\newblock {\em Discrete and Continuous Dynamical Systems - S}, 15(4):727--746.

\bibitem{pardoux1994backward}
{\'E}tienne Pardoux and Shige Peng.
\newblock Backward doubly stochastic differential equations and systems of
  quasilinear spdes.
\newblock {\em Probability Theory and Related Fields}, 98(2):209--227, 1994.

\bibitem{pardoux2005backward}
Etienne Pardoux and Shige Peng.
\newblock Backward stochastic differential equations and quasilinear parabolic
  partial differential equations.
\newblock In {\em Stochastic Partial Differential Equations and Their
  Applications: Proceedings of IFIP WG 7/1 International Conference University
  of North Carolina at Charlotte, NC June 6--8, 1991}, pages 200--217.
  Springer, 2005.

\bibitem{Peng1990}
S.~Peng.
\newblock A general stochastic maximum principle for optimal control problems.
\newblock {\em SIAM J. Control Optim.}, pages 966--979, 1990.

\bibitem{peng2010backward}
Shige Peng.
\newblock Backward stochastic differential equation, nonlinear expectation and
  their applications.
\newblock In {\em Proceedings of the International Congress of Mathematicians
  2010 (ICM 2010) (In 4 Volumes) Vol. I: Plenary Lectures and Ceremonies Vols.
  II--IV: Invited Lectures}, pages 393--432. World Scientific, 2010.

\bibitem{APF}
Michael~K. Pitt and Neil Shephard.
\newblock Filtering via simulation: auxiliary particle filters.
\newblock {\em J. Amer. Statist. Assoc.}, 94(446):590--599, 1999.

\bibitem{robbins1951stochastic}
Herbert Robbins and Sutton Monro.
\newblock A stochastic approximation method.
\newblock {\em The annals of mathematical statistics}, pages 400--407, 1951.

\bibitem{Convergence-SGLD}
I.~Sato and H.~Nakagawa.
\newblock Convergence analysis of gradient descent stochastic algorithms.
\newblock {\em Proceedings of the 31st International Conference on Machine
  Learning}, pages 982--990, 2014.

\bibitem{Convergence-SGD}
A.~Shapiro and Y.~Wardi.
\newblock Convergence analysis of gradient descent stochastic algorithms.
\newblock {\em Journal of Optimization Theory and Applications}, pages
  439--454, 1996.

\bibitem{Sny}
C.~Snyder, T.~Bengtsson, P.~Bickel, and J.~Anderson.
\newblock Obstacles to high-dimensional particle filtering.
\newblock {\em Mon. Wea. Rev.}, 136:4629--4640, 2008.

\bibitem{Tang-1998}
Shanjian Tang.
\newblock The maximum principle for partially observed optimal control of
  stochastic differential equations.
\newblock {\em SIAM J. Control Optim.}, 36(5):1596--1617, 1998.

\bibitem{Tong_EnKF}
Xin~T. Tong, Andrew~J. Majda, and David Kelly.
\newblock Nonlinear stability and ergodicity of ensemble based {K}alman
  filters.
\newblock {\em Nonlinearity}, 29(2):657--691, 2016.

\bibitem{Wang-Wu-Xiong-2018}
Guangchen Wang, Zhen Wu, and Jie Xiong.
\newblock {\em An introduction to optimal control of {FBSDE} with incomplete
  information}.
\newblock SpringerBriefs in Mathematics. Springer, Cham, 2018.

\bibitem{Yong_control}
Jiongmin Yong and Xun~Yu Zhou.
\newblock {\em Stochastic controls}, volume~43 of {\em Applications of
  Mathematics (New York)}.
\newblock Springer-Verlag, New York, 1999.
\newblock Hamiltonian systems and HJB equations.

\bibitem{zakai}
Moshe Zakai.
\newblock On the optimal filtering of diffusion processes.
\newblock {\em Z. Wahrscheinlichkeitstheorie und Verw. Gebiete}, 11:230--243,
  1969.

\bibitem{ZhangJ_BSDE}
Jianfeng Zhang.
\newblock A numerical scheme for {BSDE}s.
\newblock {\em Ann. Appl. Probab.}, 14(1):459--488, 2004.

\bibitem{Zhao_multi}
Weidong Zhao, Yu~Fu, and Tao Zhou.
\newblock New kinds of high-order multistep schemes for coupled forward
  backward stochastic differential equations.
\newblock {\em SIAM J. Sci. Comput.}, 36(4):A1731--A1751, 2014.

\bibitem{zhao2017high}
Weidong Zhao, Tao Zhou, and Tao Kong.
\newblock High order numerical schemes for second-order fbsdes with
  applications to stochastic optimal control.
\newblock {\em Communications in Computational Physics}, 21(3):808--834, 2017.

\end{thebibliography}

\end{document}